\documentclass[nohypdvips]{siamart0516}
\usepackage{amsmath}
\usepackage{amsfonts}
\usepackage{amssymb}
\usepackage{amsbsy}
\usepackage{bm}
\usepackage{graphicx}
\usepackage{extarrows}
\usepackage{mathtools}

\newcommand{\ie}{i.\,e.}%
\newcommand{\eg}{e.\,g.}%
\newcommand{\etal}{et\,al.}%
\newcommand{\wrt}{w.r.t.}%

\newcommand{\formComma}{\,\text{,}}
\newcommand{\formPeriod}{\,\text{.}}

\newcommand{\R}{\mathbb{R}}%
\newcommand{\xb}{\mathbf{x}}%
\newcommand{\nb}{\mathbf{n}}%
\newcommand{\eb}{\mathbf{e}}%
\newcommand{\gb}{\mathbf{g}}%
\newcommand{\Xb}{\mathbf{x}}%
\newcommand{\XXb}{\extendDomain{\mathbf{x}}}%

\newcommand{\alphab}{\bm{\alpha}}%
\newcommand{\betab}{\bm{\beta}}%

\newcommand{\pb}{\mathbf{p}}%
\newcommand{\p}{\textup{p}}%
\newcommand{\q}{\textup{q}}%
\newcommand{\qb}{\mathbf{q}}%
\newcommand{\vb}{\mathbf{q}}%

\newcommand{\alphav}{\underline{\bm{\alpha}}} %
\newcommand{\betav}{\underline{\bm{\beta}}} %
\newcommand{\MPDt}{\underline{\underline{\tensor{M}}}} %

\newcommand{\fv}{\underline{\mathbf{f}}}
\newcommand{\pv}{\underline{\mathbf{p}}}

\newcommand{\norm}[1]{\lVert#1\rVert}%
\newcommand{\scalarprod}[1]{\big\langle{#1}\big\rangle}%
\newcommand{\Scalarprod}[1]{\left\langle{#1}\right\rangle}%

\newcommand{\dif}{\textup{d}}
\newcommand{\exd}{\mathbf{d}} %

\newcommand{\dt}{\partial_t}

\newcommand{\dS}{\,\dif{\surf}}
\newcommand{\dV}{\,\dif{V}}

\newcommand{\fdif}{\operatorname{\delta}\!}
\newcommand{\Fdif}[2]{\frac{\fdif{#1}}{\fdif{#2}}}%

\newcommand{\FF}{\mathrm{F}}
\newcommand{\F}[1]{\FF_\mathrm{#1}}

\newcommand{\Sp}{\mathbb{S}^2}
\newcommand{\ellipsoid}{\mathcal{S}^{E}}
\newcommand{\surf}{\mathcal{S}}
\newcommand{\domain}{\Omega}

\newcommand{\Tangent}{\mathsf{T}}
\newcommand{\ProjectSurf}{\pi_{\Tangent\surf}}

\newcommand{\surfNormal}{\boldsymbol{\nu}}
\newcommand{\surfNormalI}{\nu}

\newcommand{\meanCurvature}{\mathcal{H}}
\newcommand{\gaussianCurvature}{\kappa}

\newcommand{\Grad}{\operatorname{grad}}
\newcommand{\Div}{\operatorname{div}}%
\newcommand{\Rot}{\operatorname{rot}}%

\newcommand{\laplaceBeltrami}{\Delta_{\surf}}
\newcommand{\vecLaplace}{\boldsymbol{\Delta}}
\newcommand{\laplaceDeRham}{\vecLaplace^{\textup{dR}}}
\newcommand{\laplaceDeRahm}{\laplaceDeRham}
\newcommand{\laplaceRotRot}{\vecLaplace^{\textup{RR}}}
\newcommand{\laplaceGradDiv}{\vecLaplace^{\textup{GD}}}

\newcommand{\LaplaceDeRham}{Laplace-deRham }

\newcommand{\laplaceDeRhamTilde}{\widehat{\vecLaplace}^{\textup{dR}}}
\newcommand{\NablaSurf}{\nabla_{\surf}}
\newcommand{\gDerivative}{D}

\newcommand{\laplaceDeRhamDiffuse}{\extendDomain{\vecLaplace}^{\textup{dR}}}

\newcommand{\shapeOperator}{\mathcal{B}}

\newcommand{\shapeOperatorPDt}{\underline{\underline{\shapeOperator}}} %

\newcommand{\idCoVec}{\textup{Id}_{\Tangent^{*}\surf}} %

\newcommand{\EndoTaylorLHS}{\mathcal{L}} %
\newcommand{\EndoTaylorRHS}{\mathcal{R}} %
\newcommand{\EndoTaylorLHSPDt}{\underline{\underline{\EndoTaylorLHS}}} %
\newcommand{\EndoTaylorRHSPDt}{\underline{\underline{\EndoTaylorRHS}}} %

\newcommand{\phase}{\phi}
\newcommand{\doubleWell}{W}
\newcommand{\regularization}{\zeta}
\newcommand{\doubleWellRegularized}{W_\regularization}

\newcommand{\lu}{\theta}
\newcommand{\lv}{\varphi}
\newcommand{\lxi}{\xi} %

\newcommand{\ch}[3]{\Gamma_{#1 #2}^{#3}} %
\newcommand{\lc}{E} %

\newcommand{\gext}{\extendDomain{g}}
\newcommand{\chext}[3]{\extendDomain{\Gamma}_{#1 #2}^{#3}} %
\newcommand{\lcext}{\extendDomain{\lc}}%

\newcommand{\vect}[1]{\mathbf{#1}}
\newcommand{\tensor}[1]{\mathbf{#1}}

\newcommand{\Span}[1]{\operatorname{Span}\!\left\{ #1 \right\}}

\newcommand{\landau}{\mathcal{O}} %
\newcommand{\landauxin}[1]{\landau(\lxi^{#1})}
\newcommand{\landauxi}{\landau(\lxi)}

\newcommand{\tensorfs}[3]{{{#1}_{#2}}^{#3}} %
\newcommand{\tensorsf}[3]{{{#1}^{#2}}_{#3}}

\newcommand{\K}{{K}} %
\newcommand{\Ki}{{K_1}} %
\newcommand{\Kii}{{K_2}} %
\newcommand{\Kiii}{{K_3}} %
\newcommand{\Kn}{{\omega_n}} %
\newcommand{\Kt}{{\omega_t}} %

\newcommand{\EE}{\F{\omega_n}^\surf}
\newcommand{\E}[1]{\EE[#1]}

\newcommand{\Ltwo}[2]{L^2( #1;\, #2)}%
\newcommand{\LS}{L^2(\surf)}%
\newcommand{\LtwoProd}[1]{\big({#1}\big)_{\LS}}%

\newcommand{\Csurf}[1]{{C}^{#1}(\surf)}

\newcommand{\extend}[1]{\widehat{#1}}
\newcommand{\extendDomain}[1]{\widetilde{#1}}

\newcommand{\Hdr}[2]{H^{\textup{DR}}( #1;\, #2)} %
\newcommand{\HdrExt}[2]{H^{\textup{DR}}( #1;\, #2)} %
\newcommand{\HdrDiffuse}[2]{H^{\textup{DR}}( #1;\, #2)} %

\newcommand{\Hi}[1]{H^{1}(#1)}

\newcommand{\pExt}{\extend{\pb}}

\newcommand{\qExt}{\extend{\qb}}

\newcommand{\nExt}{\extendDomain{\surfNormal}}

\newcommand{\SC}{\mathcal{K}} %
\newcommand{\Vs}{\mathcal{V}} %
\newcommand{\Es}{\mathcal{E}} %
\newcommand{\Fs}{\mathcal{T}} %
\newcommand{\face}{T} %
\newcommand{\FormSpace}{\Lambda^{1}} %
\newcommand{\PDT}{\mathfrak{T}} %

\newcommand{\numVs}{|\mathcal{V}|} %
\newcommand{\numEs}{|\mathcal{E}|} %
\newcommand{\numFs}{|\mathcal{T}|} %

\newcommand{\triangulation}{{\surf_h}}

\newcommand{\pstretch}{C}
\newcommand{\pprop}{r}
\newcommand{\ppress}{B}

\newcommand{\EuBase}[1]{\,\eb^{#1}}

\newcommand{\hidden}[1]{}

\newcommand{\shortcut}[1]{{\small\textsf{#1}}}
\newcommand{\DEC}{{\small\textsf{DEC}}}
\newcommand{\SFEM}{{\small\textsf{sFEM}}}
\newcommand{\DI}{{\small\textsf{DI}}}
\newcommand{\SPH}{{\small\textsf{SPH}}}
\usepackage[automake,nonumberlist]{glossaries}
\usepackage{glossary-mcols}

\newglossary[slg]{symbolslist}{syi}{syg}{List of symbols}

\renewcommand*{\glsgroupheading}[1]{%
\ifstrequal{A}{##1}{\relax}%
{\item[{\glsgetgrouptitle{##1}}]}%
}
\newglossarystyle{mystyle}{%
  \glossarystyle{mcolindexgroup}%
}
\setglossarystyle{mystyle}

\makeglossaries

\newglossaryentry{symb:oneform}{ name=\ensuremath{\alphab},
description={1-form, $\alphab\in\Lambda^1(\surf)$},
sort=e_symbol_dec_one-form, type=symbolslist
}
\newglossaryentry{symb:discreteoneform}{ name=\ensuremath{\alpha_h},
description={discrete 1-form, $\alpha_h\in\Lambda^1_h(\mathcal{K})$},
sort=e_symbol_dec_one-form_discrete, type=symbolslist
}
\newglossaryentry{symb:PDoneform}{ name=\ensuremath{\alphav},
description={Primal-Dual 1-form, $\alphav=(\alpha_h, \ast\alpha_h)$},
sort=e_symbol_dec_one-form_PD, type=symbolslist
}

\newglossaryentry{symb:exd}{ name=\ensuremath{\exd},
description={exterior derivative},
sort=e_symbol_derivative_exd, type=symbolslist
}
\newglossaryentry{symb:hodgestar}{ name=\ensuremath{\ast},
description={hodge star operator},
sort=e_symbol_dec_hodge_star, type=symbolslist
}

\newglossaryentry{symb:Sp}{ name=\ensuremath{\Sp},
description={unit 2-sphere},
sort=g_symbol_surface_sphere, type=symbolslist
}
\newglossaryentry{symb:ellipsoid}{ name=\ensuremath{\ellipsoid},
description={ellipsoidal surface},
sort=g_symbol_surface_ellipsoid, type=symbolslist
}
\newglossaryentry{symb:surf}{ name=\ensuremath{\surf},
description={surface, i.e., compact closed oriented Riemannian 2-dim. manifold},
sort=g_symbol_surface, type=symbolslist
}
\newglossaryentry{symb:domain}{ name=\ensuremath{\Omega},
description={domain, $\Omega\subset\R^3$},
sort=g_symbol_domain, type=symbolslist
}

\newglossaryentry{symb:projectcoords}{ name=\ensuremath{\pi},
description={coordinate projection $\pi:\Omega_\delta\to\surf$},
sort=g_symbol_projection_surface, type=symbolslist
}
\newglossaryentry{symb:projectsurf}{ name=\ensuremath{\ProjectSurf},
description={surface projection $\ProjectSurf:\Tangent\R^3\to\Tangent\surf$},
sort=g_symbol_projection_surface, type=symbolslist
}
\newglossaryentry{symb:tangentbundle}{ name=\ensuremath{\Tangent\surf},
description={tangent bundle of surface $\surf$},
sort=g_symbol_tangent_bundle, type=symbolslist
}
\newglossaryentry{symb:cotangentbundle}{ name=\ensuremath{\Tangent^{*}\surf},
description={cotangent bundle of surface $\surf$},
sort=g_symbol_tangent_co_bundle, type=symbolslist
}

\newglossaryentry{symb:surfacenormal}{ name=\ensuremath{\surfNormal},
description={outer surface normal},
sort=g_symbol_surface_normal, type=symbolslist
}

\newglossaryentry{symb:meancurvature}{ name=\ensuremath{\meanCurvature},
description={mean curvature $\meanCurvature=\Div\boldsymbol{\nu}$},
sort=g_symbol_curvature_mean, type=symbolslist
}
\newglossaryentry{symb:gaussiancurvature}{ name=\ensuremath{\gaussianCurvature},
description={Gaussian curvature},
sort=g_symbol_curvature_gaussian, type=symbolslist
}
\newglossaryentry{symb:characteristic}{ name=\ensuremath{\chi( \surf )},
description={characteristic of the surface $\surf$},
sort=g_symbol_surface_characteristic, type=symbolslist
}

\newglossaryentry{symb:grad}{ name=\ensuremath{\Grad},
description={surface gradient},
sort=d_symbol_fot_grad, type=symbolslist
}
\newglossaryentry{symb:div}{ name=\ensuremath{\Div},
description={surface divergence},
sort=d_symbol_fot_div, type=symbolslist
}
\newglossaryentry{symb:rot}{ name=\ensuremath{\Rot},
description={surface curl},
sort=d_symbol_fot_rot, type=symbolslist
}
\newglossaryentry{symb:laplacebeltrami}{ name=\ensuremath{\laplaceBeltrami},
description={surface Laplace-Beltrami operator},
sort=d_symbol_laplace_beltrami, type=symbolslist
}
\newglossaryentry{symb:laplacederham}{ name=\ensuremath{\laplaceDeRham},
description={surface Laplace-deRham operator},
sort=d_symbol_laplace_derham, type=symbolslist
}

\newglossaryentry{symb:shapeoperator}{ name=\ensuremath{\shapeOperator},
description={shape operator $\shapeOperator=-\Grad\surfNormal$},
sort=g_symbol_shape_operator, type=symbolslist
}

\newglossaryentry{symb:phase}{ name=\ensuremath{\phase},
description={phase-field variable},
sort=p_symbol_phase, type=symbolslist
}
\newglossaryentry{symb:signeddistance}{ name=\ensuremath{d_\surf(\mathbf{x})},
description={signed-distance function},
sort=p_symbol_signed_distance, type=symbolslist
}
\newglossaryentry{symb:surfacedelta}{ name=\ensuremath{\delta_{\surf}},
description={surface delta-function},
sort=p_symbol_phase_delta, type=symbolslist
}
\newglossaryentry{symb:doublewell}{ name=\ensuremath{\doubleWell},
description={double-well, $\doubleWell(\phase)\simeq\delta_\surf$},
sort=p_symbol_phase_doublewell, type=symbolslist
}
\newglossaryentry{symb:epsilon}{ name=\ensuremath{\varepsilon},
description={interface thickness of phase-field},
sort=p_symbol_phase_epsilon, type=symbolslist
}
\newglossaryentry{symb:doublewellregularization}{ name=\ensuremath{\regularization},
description={double-well regularization},
sort=p_symbol_phase_doublewell_regularization, type=symbolslist
}

\newglossaryentry{symb:lu}{ name=\ensuremath{\lu},
description={co-latitude coordinate, $\lu\in[0,\pi]$},
sort=g_symbol_coordinate_u, type=symbolslist
}
\newglossaryentry{symb:lv}{ name=\ensuremath{\lv},
description={azimuthal coordinate, $\lv\in[0,2\pi)$},
sort=g_symbol_coordinate_v, type=symbolslist
}
\newglossaryentry{symb:lxi}{ name=\ensuremath{\lxi},
description={coordinate in normal direction of the surface},
sort=g_symbol_coordinate_xi, type=symbolslist
}

\newglossaryentry{symb:christoffel}{ name=\ensuremath{\ch{i}{j}{k}},
description={Christoffel symbols of second kind},
sort=g_symbol_christoffel, type=symbolslist
}

\newglossaryentry{symb:levicivita}{ name=\ensuremath{\lc_{IJK}},
description={Levi-Civita symbols},
sort=g_symbol_levi_civita, type=symbolslist
}

\newglossaryentry{symb:metric}{ name=\ensuremath{\tensor{g}},
description={Riemannian metric tensor},
sort=g_symbol_metric_g, type=symbolslist
}
\newglossaryentry{symb:metric_det}{ name=\ensuremath{|\tensor{g}|},
description={determinant of $\tensor{g}$},
sort=g_symbol_metric_g_det, type=symbolslist
}

\newglossaryentry{symb:K}{ name=\ensuremath{\K},
description={uniform Frank constant},
sort=m_symbol_constant_frank, type=symbolslist
}
\newglossaryentry{symb:Kn}{ name=\ensuremath{\Kn},
description={penalty constant for normality},
sort=m_symbol_constant_normal, type=symbolslist
}
\newglossaryentry{symb:Kt}{ name=\ensuremath{\Kt},
description={penalty constant for tangentiality},
sort=m_symbol_constant_tangential, type=symbolslist
}

\newglossaryentry{symb:defectfusion}{ name=\ensuremath{\epsilon_f},
description={error in the defect fusion time},
sort=m_symbol_error_defect_fusion_time, type=symbolslist
}
\newglossaryentry{symb:energyerror}{ name=\ensuremath{\epsilon_e},
description={(normalized) mean energy error},
sort=m_symbol_error_energy, type=symbolslist
}

\newglossaryentry{symb:EE}{ name=\ensuremath{\EE},
description={weak surface Frank-Oseen energy},
sort=m_symbol_energy_weakFrank, type=symbolslist
}

\newglossaryentry{symb:SC}{ name=\ensuremath{\SC},
description={simplicial complex},
sort=e_symbol_dec_simplicial_complex, type=symbolslist
}
\newglossaryentry{symb:Vs}{ name=\ensuremath{\Vs},
description={set of vertices, with number $\numVs$},
sort=e_symbol_dec_vertices, type=symbolslist
}
\newglossaryentry{symb:Es}{ name=\ensuremath{\Es},
description={set of edges, with number $\numEs$},
sort=e_symbol_dec_edges, type=symbolslist
}
\newglossaryentry{symb:Fs}{ name=\ensuremath{\Fs},
description={set of faces, with number $\numFs$},
sort=e_symbol_dec_faces, type=symbolslist
}

\newglossaryentry{symb:edge}{ name=\ensuremath{e},
description={edge, $e\in\Es$},
sort=e_symbol_dec_edge, type=symbolslist
}
\newglossaryentry{symb:dualedge}{ name=\ensuremath{\star e},
description={dual edge of $e$ (voronoi edge)},
sort=e_symbol_dec_edge_dual, type=symbolslist
}
\newglossaryentry{symb:edgevec}{ name=\ensuremath{\mathbf{e}},
description={edge vector along edge $e$},
sort=e_symbol_dec_edgevec, type=symbolslist
}
\newglossaryentry{symb:dualedgevec}{ name=\ensuremath{\mathbf{e}_\star},
description={dual edge vector along dual chain $\star e$},
sort=e_symbol_dec_edgevec_dual, type=symbolslist
}

\newglossaryentry{symb:vertex}{ name=\ensuremath{v},
description={vertex, $v\in\Vs$},
sort=e_symbol_dec_vertex, type=symbolslist
}
\newglossaryentry{symb:dualvertex}{ name=\ensuremath{\star v},
description={dual vertex (voronoi cell)},
sort=e_symbol_dec_vertex_dual, type=symbolslist
}

\newglossaryentry{symb:face}{ name=\ensuremath{\face},
description={face, $\face\in\Fs$},
sort=e_symbol_dec_face, type=symbolslist
}

\newglossaryentry{symb:flat}{ name=\ensuremath{\flat},
description={lowering indices},
sort=e_symbol_dec_musical_flat, type=symbolslist
}
\newglossaryentry{symb:sharp}{ name=\ensuremath{\sharp},
description={rising indices},
sort=e_symbol_dec_musical_sharp, type=symbolslist
}

\newglossaryentry{symb:time}{ name=\ensuremath{t_k},
description={discrete time},
sort=m_symbol_time, type=symbolslist
}
\newglossaryentry{symb:tau}{ name=\ensuremath{\tau_k},
description={time step width in th $k$-th time step},
sort=m_symbol_time_step, type=symbolslist
}

\makeatother

\title{Orientational order on surfaces -- the coupling of topology, geometry, and dynamics}

\author{M. Nestler\footnotemark[2]\ \footnotemark[3]
\and I. Nitschke\footnotemark[3]
\and S. Praetorius\footnotemark[3]
\and A. Voigt\footnotemark[3]}

\begin{document}
\maketitle

\renewcommand{\thefootnote}{\fnsymbol{footnote}}
\footnotetext[2]{Corresponding author: michael.nestler@tu-dresden.de (Michael Nestler)}
\footnotetext[3]{Institut f\"{u}r Wissenschaftliches Rechnen, Technische Universit\"{a}t Dresden, Zellescher Weg 12--14, 01062 Dresden, Germany (michael.nestler@tu-dresden.de, ingo.nitschke@tu-dresden.de, simon.praetorius@tu-dresden.de, axel.voigt@tu-dresden.de)}
\renewcommand{\thefootnote}{\arabic{footnote}}

\begin{abstract}
We consider the numerical investigation of surface bound orientational order using unit tangential vector fields by means of
a gradient-flow equation of a weak surface Frank-Oseen energy. The energy is composed of intrinsic and extrinsic contributions,
as well as a penalization term to enforce the unity of the vector field.
Four different numerical discretizations, namely a discrete exterior calculus approach, a method based on vector spherical
harmonics, a surface finite-element method, and an approach utilizing an implicit surface description, the diffuse interface
method, are described and compared with each other for surfaces with Euler characteristic 2. We demonstrate the influence
of geometric properties on realizations of the Poincar\'e-Hopf theorem and show examples where the energy is decreased by
introducing additional orientational defects.

\end{abstract}



\pagestyle{myheadings}
\thispagestyle{plain}
\markboth{M. Nestler \etal}{Orientational order on surfaces}

\section{Introduction}\label{sec:Introduction}
We consider surface bound systems of densely packed rod like particles that tend to align tangentially. The systems are modeled by a mesoscopic field theoretical description using an average direction and an order parameter, measuring the local variance of alignment towards this average direction. In flat space an uniformly ordered ground state can be established. This is no longer true for curved space, which induces distortions of this ground state, eventually inhibiting the propagation of preferred orientational order throughout the whole system. This leads to the emergence of defects, which for surfaces $\surf$ with Euler characteristic $\chi(\surf) \neq 0$ is a consequence of the Poincar\'e-Hopf theorem. However, the type of the defects, their number, as well as their position are mostly unknown. The realization of the Poincar\'e-Hopf theorem depends on geometric properties of the surface and dynamics of the evolution. It is the goal of this paper to provide numerical methods to explore these interesting and nontrivial connections between topology, geometry and dynamics. Besides the mathematical issues, the problem is of interest in the physics and materials science community due to its envisioned technological applications \cite{Nelson2002}.

We focus on orientational ordering in polar order dynamics. The model follows as limit of a thin film formulation of a modified Frank-Oseen energy \cite{Frank1958} and is formulated as an $L^2$-gradient flow, which leads to a vector-valued partial differential equation on the surface. Previous work has postulated a purely intrinsic formulation, extending the flat space model to curved space \cite{Nelson1983,Lubensky1992,Leon2011}.  More recent research \cite{Napoli2012a,Napoli2012b,Segatti2014} derives a surface Frank-Oseen energy as limit of a thin film formulation. This approach adds to the intrinsic model an explicit influence of the embedding space by extrinsic quantities. However, the limit is only established for surfaces with $\chi(\surf) = 0$ and only allows defect free configurations. All approaches focus only on the steady state and utilize continuous optimization methods \cite{Kralj2011} or Monte-Carlo based methods \cite{selinger2011,li2014,nguyen2013} to evaluate the minimizers. To complement these models and methods we derive a more general thin film limit, valid also for surfaces with $\chi(\surf) \neq 0$ and focus on the dynamics of orientational order on such surfaces.

Starting from the general surface modeling provided in \autoref{sec:ModelDerivation} we establish suitable reformulations to apply different numerical methods and solve the resulting dynamic equations. We propose methods based on a coordinate free
framework as well as methods adapted for the Cartesian coordinates of the embedding $\R^3$ by using a penalty term approach. \autoref{sec:Notation} gives the general notations and \autoref{sec:NumericalMethods} presents the methods of discrete exterior calculus (\DEC), vector spherical harmonics (\SPH),
surface finite elements (\SFEM) and diffuse interface modeling (\DI). We compare results of these methods in \autoref{sec:ComputationalResults} to provide estimations on numerical quality and computational cost. Further, we use these methods to perform experiments investigating the influence of geometry on emergence and energetical stability of non-minimal defect configurations and demonstrate the possibility to decrease the energy by introducing additional defects. The model formulations and proposed
methods will provide a modeling and numerical toolkit ready to be applied to polar orientational order in curved space and related physical systems out of equilibrium. This and the implication for solving vector-valued partial differential equations on surfaces will be discussed in \autoref{sec:ConclusionOutlook}.

\section{Model derivation}\label{sec:ModelDerivation}
  Two major continuous theories to describe orientational order in liquid crystals exist. On the one hand, the Frank-Oseen theory uses a vector field to describe average molecular ordering,
  while, on the other hand, the Landau-de Gennes theory is based on a matrix expression (called Q-tensor). Both models are widely used and indeed coincide in flat 2D space for a specific
  set of elastic terms, see \cite{ball2011, iyer2015}. Besides this agreement, the Frank-Oseen modeling can not account for a physical head-to-tail symmetry of the material, which is naturally considered in the Landau-de Gennes theory. For a mathematical review on both modeling approaches we refer to \cite{Ball2016}. Due to its relative simplicity
  we here consider only the Frank-Ossen theory as a modeling framework. Being aware of the fact that additional physical effects will occur within a corresponding Landau-de Gennes theory.

  In our framework the average alignment of anisometric molecules can be expressed by a unit vector  \( \pb \), in the following called director, that represents the direction of the average alignment axis. In order to describe the spatial variation of a director field a free energy $\FF$ can be formulated that incorporates energy costs due to spatial distorsions. The energy reads in simplified form \cite{Oswald2005}
  \begin{align}
    \F{F}\left[ \pb, \Omega \right] &= \frac{1}{2}\int_{\Omega} \Ki\left( \nabla\cdot\pb \right)^{2}
                        + \Kii\left( \pb\cdot\left( \nabla\times\pb \right) \right)^{2}
                        + \Kiii\left\| \pb\times\left( \nabla\times\pb \right) \right\|^{2} \dV\,,
  \end{align}
  with $K_1, K_2$, and $K_3$ the Frank phenomenological constants and \( \Omega \subset \R^{3} \) a three dimensional domain. The functional $\F{F}$ contains three contributions related to
  deformations of $\pb$, namely (from left to right) for splay, twist, and bend. We here consider the one-constant approximation \( \K:=\Ki=\Kii=\Kiii \). The distortion energy thus reads
  \begin{align}
    \F{OC}\left[ \pb, \Omega \right] &= \frac{\K}{2}\int_{\Omega} \left( \nabla\cdot\pb \right)^{2} + \| \nabla \times \pb \|^2 \dV\,.
  \end{align}
  To arrive at a surface formulation, we consider a thin shell \( \Omega = \Omega_{\delta} \) around a compact smooth Riemannian surface $\surf$, with thickness $\delta$ sufficiently small, and $\pb$ parallel to the surface and parallel transported in normal direction to the surface. The limiting case of $\F{OC}\left[ \pb, \Omega_\delta \right]$,
  $\delta \searrow 0$, where $\Omega_{\delta}$ collapses to the surface, has been considered in \cite{Napoli2012b} for surfaces with \( \chi( \surf ) = 0 \) and thus only for defect
  free configurations. This result cannot simply be extended to more general surfaces, as a smooth vector field with unit norm exists if and only if \( \chi( \surf ) = 0 \). This topological result can also be extended to the corresponding Sobolev space \cite{Segatti2016} and thus turns out to be useless for any investigation of defects in unit vector fields on surfaces. While in mathematical terms these defects can be considered as discontinuities, in physical terms the liquid crystal undergoes a phase transition to an isotropic phase at the defect. To enable a continuous director field $\pb$ and to incorporate this phase transition, we drop the constraint $\| \pb \| = 1 $ and consider $\|\pb\|$ as an order parameter. This parameter ranges from $0$, describing the isotropic phase, to $1$, for the ordered phase of the liquid crystal. To enforce a prevalent ordered phase, we add a well-known quartic state potential to the free energy with penalty parameter $\Kn$. It is evident that the radius of the defect core, the domain where the local alignment breaks down, is closely connected to $\Kn$. Since we are interested in orientational ordering of a prevalent ordered state, we choose $\Kn \gg \K$, effectively enforcing defects with small core radius. The corresponding energy reads
    \begin{align}
    \label{eq:FOEPenalized}
    \F{\Kn}\left[ \pb, \Omega \right] &= \frac{\K}{2}\int_{\Omega} \left( \nabla\cdot\pb \right)^{2} + \| \nabla \times \pb \|^2 \dV + \frac{\Kn}{4}\int_{\Omega} \left( \| \pb \|^2 -1 \right)^{2} \dV\,.
  \end{align}
 Extending the ansatz of \cite{Napoli2012b}, the limit $\delta \searrow 0$ can now be considered also for \( \chi( \surf ) \neq 0 \), see \autoref{sec:Limit}. We obtain $\lim_{\delta\searrow 0} \frac{1}{\delta}\F{\Kn}\left[ \pb,
 \Omega_{\delta} \right] =\E{\pb}$, which we call the \textit{weak surface Frank-Oseen energy}
 \begin{align}
    \E{\pb} &= \F{I}^{\surf}[\pb] +  \F{E}^{\surf}[\pb] + \frac{\Kn}{4}\int_{\surf} \left( \| \pb \|^2 -1 \right)^{2} \dS.
  \end{align}
It consists of an intrinsic contribution $\F{I}^{\surf}[\pb]$ and an extrinsic contribution $\F{E}^{\surf}[\pb]$ to the distortion energy, as in
  \cite{Napoli2012a,Napoli2012b}, and the additional penalty term, which contains the 2-norm $\|\cdot\|$.
  In the following we assume $\pb\in \Tangent\surf$ the tangent bundle of $\surf$. Then, the intrinsic distortion energy $\F{I}^\surf$ can be expressed in terms of the surface divergence ``$\Div$'' and the surface curl ``$\Rot$'' of $\pb$:
  \begin{align}
    \F{I}^\surf[\pb] &= \frac{\K}{2}\int_{\surf} \left( \Div\pb \right)^{2} + \left( \Rot\pb \right)^{2} \dS\, .
  \end{align}
  Introducing further the shape operator $\shapeOperator=-\Grad\surfNormal$ of $\surf$ with outer surface normal $\surfNormal$, the extrinsic contributions can be written as
  \begin{align}
   \F{E}^{\surf}[\pb] &= \frac{\K}{2}\int_{\surf} \| \shapeOperator \cdot \pb \|^{2}  \dS.
  \end{align}
Putting all parts together, we finally obtain
  \begin{align}\label{eq:SurfEnergy}
    \E{\pb} &= \frac{\K}{2}\int_{\surf} \left( \Div\pb \right)^{2} + \left( \Rot\pb \right)^{2} +  \| \shapeOperator \cdot \pb \|^{2} \dS
                                +\frac{\Kn}{4}\int_{\surf} \big( \left\| \pb \right\|^{2} - 1 \big)^{2} \dS.
  \end{align}

  For the description of the minimization of $\E{\pb}$, we define the function spaces
  \begin{align*}
    H(\Div,\surf, \Tangent \surf) &:= \big\{\pb \in \Ltwo{\surf}{\Tangent\surf}\,:\;\Div\pb\in\LS\big\}\formComma \\
    H(\Rot,\surf, \Tangent \surf) &:= \big\{\pb \in \Ltwo{\surf}{\Tangent\surf}\,:\;\Rot\pb\in\LS\big\}\formComma
  \end{align*}
  and furthermore the space $\Hdr{\surf}{\Tangent \surf} := H(\Div,\surf, \Tangent \surf)\cap H(\Rot,\surf, \Tangent\surf)$. The minimization of the weak surface Frank-Oseen energy reads
  \begin{align*}
    \pb^\ast &= \operatorname{argmin} \big\{ \E{\pb}\,:\; \pb\in \Hdr{\surf}{\Tangent \surf}\big\}\formPeriod
  \end{align*}
  In \cite{Chen1989,Segatti2014} the convergence of minimizers of $\EE$ to the sharp energy $\F{I}^\surf[\pb] + \F{E}^\surf[\pb]$, as $\Kn\rightarrow\infty$, is
  analyzed and proven for the case \( \chi( \surf ) = 0 \).

  Dynamical equations to minimize the functional $\EE$ can be formulated by
  means of an \( L^{2} \)-gradient flow approach,
  \[
    \dt\pb = - \Fdif{\EE}{\pb}[\pb]\,,
  \]
  where the gradient of $\EE$ has to be interpreted \wrt\ the $\Ltwo{\surf}{\Tangent\surf}$-inner product.
  For $\vb\in \Hdr{\surf}{\Tangent \surf}$ this reads
  \begin{align*}
    \MoveEqLeft[3]{\int_{\surf} \Scalarprod{\Fdif{\EE}{\pb}[\pb],\vb} \dS} & \\
	&= \int_{\surf} -\K\left(\Div\pb\Div\vb + \Rot\pb\Rot\vb\right)
        +\K\scalarprod{\shapeOperator\pb,\shapeOperator\vb}
	      + \Kn \big( \left\| \pb \right\|^{2} - 1 \big) \scalarprod{\pb,\vb} \dS \\
	&= \int_{\surf} \K\scalarprod{\laplaceDeRham\pb,\vb}
        +\K\scalarprod{\shapeOperator^{2}\pb,\vb}
	      + \Kn \big( \left\| \pb \right\|^{2} - 1 \big) \scalarprod{\pb,\vb} \dS\formComma
  \end{align*}
  with $\laplaceDeRham$ the \LaplaceDeRham operator.
  This leads to the evolution equation
  \begin{align}\label{eq:VectorPDE}
    \partial_{t}\pb + \K\left(\laplaceDeRham\pb + \shapeOperator^{2}\pb\right)
        + \Kn \left( \left\| \pb \right\|^{2} - 1 \right)\pb &= 0 \formComma\quad\text{ in }\surf\times\left( 0,\infty \right)
  \end{align}
  with the initial condition \( \pb\left( t=0  \right) = \pb^{0} \in \Tangent\surf \).
  The gradient flow approach guarantees dissipative dynamics and stationary solutions of
  \eqref{eq:VectorPDE} as local minima of $\EE$. Note that the
  sign of the vectorial Laplacian is different from the sign of the scalar Laplacian found
  in classical diffusion-like equations, since we follow the convention of \cite{Marsden1988}.

  Introducing the covariant director \( \alphab := \pb^{\flat} \in \Tangent^{*}\surf\), an equivalent formulation
  of equation \eqref{eq:VectorPDE} in terms of its dual vectors can be stated:
  \begin{align}\label{eq:FormPDE}
    \partial_{t}\alphab + \K\left(\laplaceDeRham\alphab + \shapeOperator^{2}\alphab\right)
        + \Kn \left( \left\| \alphab \right\|^{2} - 1 \right)\alphab &= 0 \formComma
  \end{align}
  with \( \alphab^{0} = (\pb^{0})^{\flat}\in \Tangent^{*}\surf  \), where we have used the notation of a
  musical isomorphism $\flat$ to denote the flattening operation. Both formulations of the gradient-flow problem, \eqref{eq:VectorPDE} and \eqref{eq:FormPDE}, are
  implemented in the present paper by means of several numerical approaches.

\section{Notation}\label{sec:Notation}
We consider a compact closed oriented Riemannian 2-dimensional manifold $\surf\subset\R^3$ parametrized by
the local coordinates $ \lu, \lv$:
\begin{equation}\label{eq:parametrization}
  \Xb:\R^{2}\supset U \rightarrow \R^{3};\ \left( \lu, \lv \right)\mapsto \Xb \left( \lu, \lv \right)\, .
\end{equation}
Thus the embedded \( \R^{3} \) representation of the surface is given by \( \surf = \Xb(U) \). The unit
outer normal of $\surf$ at point $\Xb$ is denoted by $\surfNormal(\Xb)$. An implicit description of the
surface is given by the signed-distance function
\begin{equation}\label{eq:implicit-representation}
d_\surf(\extendDomain{\Xb}) := \left\{\begin{array}{ll}
  -\inf_{\mathbf{y}\in\surf} \|\extendDomain{\Xb} - \mathbf{y}\| & \text{ for } \extendDomain{\Xb}\in G  \\
  \inf_{\mathbf{y}\in\surf} \|\extendDomain{\Xb} - \mathbf{y}\| & \text{ for } \extendDomain{\Xb}\in \R^{3}\setminus \bar{G}
  \end{array}\right.\formComma
\end{equation}
with a bounded open set \( G\subset\R^{3} \) and \( \partial G = \surf \).
The corresponding extended surface normal $\extendDomain{\surfNormal}:\R^3\to\R^3$ can be
calculated by
\begin{equation}
\extendDomain{\surfNormal} := \frac{\nabla d_\surf}{\|\nabla d_\surf\|},\quad\text{ with }\extendDomain{\surfNormal}\big|_\surf = \surfNormal\text{ and }\|\nabla d_\surf\|=1\formComma
\end{equation}
see, \eg, \cite{Dziuk2013}.

The key ingredient in differential geometry and tensor analysis on Riemannian manifolds is
the positive definite metric tensor
\begin{align}
  \tensor{g} =
    \begin{bmatrix}
      g_{\lu\lu} & g_{\lu\lv}\\
      g_{\lu\lv} & g_{\lv\lv}
    \end{bmatrix}
    = g_{\lu\lu}\,d\lu^{2} + 2g_{\lu\lv}\,d\lu \,d\lv + g_{\lv\lv}\,d\lv^{2} \formPeriod
\end{align}
The covariant components of the metric tensor are given by \( \R^{3} \) inner products of partial derivatives of \( \Xb \),
\ie, \( g_{ij} = \partial_{i}\Xb\cdot\partial_{j}\Xb \).
The components of the inverse tensor \( \tensor{g}^{-1} \) are denoted by \( g^{ij} \) and the determinant of
\( \tensor{g} \) by \( \left| \tensor{g} \right| \).
We denote by \( \left\{ \partial_{\lu}\Xb,\partial_{\lv}\Xb  \right\}\) the canonical basis to describe contravariant (tangential) vectors
\( \pb(\Xb)\in \Tangent_{\Xb}\surf \), \ie, \( \pb = p^{\lu} \partial_{\lu}\Xb + p^{\lv}\partial_{\lv}\Xb \) at a point \( \Xb\in\surf \).
Furthermore, with the arising dual basis \( \left\{ d\lu, d\lv \right\} \) we are able to write an arbitrary 1-form (covariant
vector) \( \alphab\in \Tangent^{*}_{\Xb}\surf \) as \( \alphab = p_{\lu}d\lu + p_{\lv}d\lv \).
This identifier choice of the covariant vector coordinates \( p_{i} \) in conjunction with representation of \( \pb
\) as above implies that \( \alphab \) and \( \pb \) are related by \( \alphab = \pb^{\flat} \) and \( \pb = \alphab^{\sharp} \), respectively.
Explicitly lowering and rising the indices can be done using the metric tensor \( \tensor{g} \) by
\( p_{i} = g_{ij}p^{j} \) and \( p^{i} = g^{ij}p_{j} \), respectively.

In a (tubular) neighborhood $\Omega_\delta$ of $\surf$, defined by $\Omega_\delta := \{ \extendDomain{\Xb}\in\R^3\,:\,d_\surf(\extendDomain{\Xb}) < \frac{1}{2}\delta \}$,
a coordinate projection $\Xb\in\surf$ of $\extendDomain{\Xb}\in\R^3$ is introduced, such that
\begin{equation}
\extendDomain{\Xb} = \Xb + d_\surf(\extendDomain{\Xb})\surfNormal(\Xb).
\end{equation}
For $\delta$ sufficiently small (depending on the local curvature of the surface),
this projection is injective, see \cite{Dziuk2013}. For a given $\extendDomain{\Xb}\in\Omega_\delta$ the coordinate
projection of $\extendDomain{\Xb}$ will also be called gluing map, denoted by $\pi:\Omega_\delta\to\surf,\,\extendDomain{\Xb}\mapsto\Xb$.

Scalar functions $f:\surf\to\R$ and vector fields $\pb:\surf\to\Tangent\surf$ can be smoothly extended
in the neighborhood $\Omega_\delta$ of $\surf$ by utilizing the coordinate projection, \ie,
extended fields $\tilde{f}:\Omega_\delta\to\R$ and $\tilde{\pb}:\Omega_\delta\to\R^3$ are defined by
\begin{equation}\label{eq:smooth_extension}
\tilde{f}(\extendDomain{\Xb}) := f(\Xb)\quad\text{ and }\quad\tilde{\pb}(\extendDomain{\Xb}) := \pb(\Xb)\formComma
\end{equation}
respectively, for $\extendDomain{\Xb}\in\Omega_\delta$ and $\Xb$ the corresponding coordinate projection. This extension can be realized by implementing a Hopf-Lax
formula on discrete grids representing the surface and its neighborhood, similar to a redistancing method,
see \cite{Bornemann2006, Burger2008}.

\subsection{Function spaces}
For scalar fields $f,g:\surf\to\mathbb{K}\in\{\R, \mathbb{C}\}$ and vector fields $\pb,\qb:\surf\to\Tangent\surf$ an $L^2$ inner product is given by
\begin{align}
\left(f,\,g\right)_{\LS} &:= \int_\surf f\,\bar{g}\dS\formComma \\
\left(\pb,\,\qb\right)_{\Ltwo{\surf}{\Tangent\surf}} &:= \int_\surf \Scalarprod{\pb,\,\bar{\qb}}\dS\formComma
\end{align}
respectively, with $\bar{g},\bar{\qb}$ the complex conjugates\footnote{In the spherical harmonics method the functions are complex-valued and thus, we need a complex $L^2$ inner product. For all real-valued functions the complex conjugation can be ignored.} and $\Scalarprod{\cdot,\,\cdot}$ the local inner product, see \autoref{tabFirstOrder}.
These $L^2$ inner products define the corresponding $\LS$ and $\Ltwo{\surf}{\Tangent\surf}$ Hilbert spaces, respectively.

\subsection{Differential calculus}
There are many ways to describe classical differential operators on surfaces.
The choice of representation arises from the context that we want to use.
In \autoref{tabFirstOrder} first order differential operations on scalars and vector fields
and an inner product are summarized and listed for the specific context.
\begin{table}\centering
\renewcommand{\arraystretch}{1.5}
\begin{tabular}{|c|c|c|c|}\hline
  Symbolic & Local coord. & \( \R^{3} \) coord. & EC\\\hline
  \( \left\langle \pb, \qb \right\rangle \)
          & \( p_{i}q^{i} \)
              & \( \pExt\cdot\qb \)
      & \( *\left( \alphab \wedge *\betab \right) \)
              \\\hline %
  \( \Grad f \)
          & \( g^{ij}\partial_{j}f\partial_{i}\Xb \)
      & \( \ProjectSurf\nabla f \)
      & \(  \exd f  \)
      \\\hline
  \( \Rot f \)
          & \( \frac{1}{\sqrt{|\tensor{g}|}} \left( \partial_{\lu} f \partial_{\lv}\Xb -\partial_{\lv} f \partial_{\lu}\Xb\right) \)
      & \( \surfNormal \times \nabla f \)
      & \(  *\exd f  \)
      \\\hline
  \(  \Div \pb \)
          & \( \partial_{i}p^{i} + \frac{1}{\sqrt{|\tensor{g}|}}p^{i}\partial_{i}\sqrt{|\tensor{g}|} \)
      & \( \nabla \cdot \pExt - \surfNormal \cdot (\nabla \pExt \cdot \surfNormal)\)
      & \( *\exd*\alphab \)
      \\\hline
  \( \Rot \pb \)
          & \( \frac{1}{\sqrt{|\tensor{g}|}} \left( \partial_{\lu}  p_{\lv} -\partial_{\lv} p_{\lu}\right)  \)
      & \( \left( \nabla\times\pExt \right)\cdot\surfNormal \)
      & \( *\exd \alphab \)
      \\\hline
\end{tabular}
\caption{Various representations of the inner product and first order differential operators on surfaces for scalar fields
$f:\surf\rightarrow\mathbb{R}$ and tangential vector fields $\pb,\qb:\surf\rightarrow\Tangent\surf$
        or \( \R^{3} \) vector fields \( \pExt:\surf\rightarrow\Tangent\R^{3}\cong \R^{3} \) are listed.
Vector-valued images are represented in a contravariant form.
In the formulation in $\mathbb{R}^3$ coordinates the scalar field
$f$ and vector field $\pExt$ with respect to the Euclidean basis \( \{\EuBase{x},\EuBase{y},\EuBase{z}\} \)
are assumed to be defined in a neighborhood of $\surf$.
In the column ``Exterior Calculus'' (EC) all is in the space of 1-forms that are related to the vector fields $\pb$ by  \( \alphab =
\pb^{\flat} \), \( \betab = \qb^{\flat} \) and the images can be compared with other columns by rising the indices.
}
\label{tabFirstOrder}
\renewcommand{\arraystretch}{1}
\end{table}
With introduced local coordinate chart above, we can use the inner metric \( \tensor{g} \) and partial derivatives \( \partial \)
(column ``Local coord.'' in \autoref{tabFirstOrder}).
In the Euclidean space \( \R^{3} \), where the surface is embedded, it is possible to describe the differential operators using
\( \R^{3} \) operators like \( \cdot \), \( \times \) or \( \nabla \) and the surface normal \( \surfNormal \).
The extension from the surface \( \surf \) to \( \R^{3} \) rises some choices of
embedding the \( \R^{3} \) vector space structure to the tangential bundle of the surface.
We use in this paper a pointwise defined normal projection
\begin{align}\label{eq:surfaceProjection}
\begin{aligned}
  \ProjectSurf(\xb):  \Tangent_\xb\R^3 \cong \R^{3} &\rightarrow \Tangent_\xb\surf;\\
  \pExt(\xb) &\mapsto \pExt(\xb) - \surfNormal(\xb)(\surfNormal(\xb)\cdot \pExt(\xb)) = \pb(\xb)
\end{aligned}
\end{align}
for all \( \xb \in \surf \), which maps an \( \R^{3} \) vector
\( \pExt = p_{x}\EuBase{x} + p_{y}\EuBase{y} + p_{z}\EuBase{z} \in \R^3  \), not necessarily tangential to the surface,
to a tangential vector \( \pb \in \Tangent_\xb\surf \). We drop the argument $\xb$ when applied to vector fields living on $\surf$.
Some flexibility arises in the choice of the first order differential operators for non-tangential vector fields defined on $\surf$,
see the operators listed in column ``\( \R^{3} \) coord.'' in \autoref{tabFirstOrder}. With this notation we can express the shape
operator as a linear map $\shapeOperator=\{\shapeOperator_{\, j}^i\}:\Tangent\surf\rightarrow\Tangent\surf$ in local and $\R^3$ coordinates by
\begin{align}
 \shapeOperator_{\, j}^i = -g^{ik}\left(\partial_{j}\surfNormal\cdot \partial_{k}\Xb \right) \quad i,j,k = 1,2 \quad \mbox{and} \quad \shapeOperator_{i j} = -\left[ \Grad \surfNormalI_j \right]_i \quad i,j = 1,2,3 \formComma
\end{align}
respectively.
This operator is symmetric, \ie, \( \left\langle \qb, \shapeOperator \pb \right\rangle = \left\langle \pb, \shapeOperator \qb \right\rangle \)
for all \( \pb,\qb\in\Tangent\surf \).
For the shape operator on the dual space in local coordinates
\begin{align}
   ^{\flat}\shapeOperator^{\sharp} = \{g_{ik}\shapeOperator_{\, l}^k g^{lj}\} = \{\shapeOperator_{i}^{\,\, j}\}
      :\Tangent^{*}\surf\rightarrow\Tangent^{*}\surf\formComma
\end{align}
we will omit the superscripts \( \sharp \) and \( \flat \) and write \( \shapeOperator \) shortly,
if it is clear on which object the shape operator is acting.
Throughout these definitions, we require the operators to coincide with surface operators for tangential fields.

From a physical point of view, neither \( \pb\in \Tangent\surf \)
nor the differential operator listed in column ``Symbolic'' in \autoref{tabFirstOrder} need explicitly defined coordinate charts.
Such a coordinate-free formulation ensures conformance in every smooth coordinate system.
In the context of exterior calculus (EC) a graded associative algebra referring to the wedge product \( \wedge \) and differential forms is introduced to implement such a coordinate-free formulation.
All fundamental first order differential operators listed in column ``EC'' in \autoref{tabFirstOrder}
can be described by the Hodge star \( * \) and the exterior derivative \( \exd \),
which arise algebraically, see \cite{Marsden1988} for details.

The Laplace operators in this paper can be obtained by composing first order operators.
The Rot-Rot-Laplace and Grad-Div-Laplace for vector-valued functions (and 1-forms) are defined by
\begin{align}
  \laplaceRotRot \pb &:= \Rot\Rot\pb
  &\text{and}&&
  \laplaceGradDiv \pb &:= \Grad\Div\pb \formPeriod
\end{align}
In \cite{Marsden1988} the \LaplaceDeRham operator is defined for \( k \)-forms on an \( n \)-dimensional Riemannian manifold by
\( \laplaceDeRham:= \left( -1 \right)^{nk+1}\left( *\exd*\exd  + \exd*\exd*\right)\).
For vector fields, we define the \LaplaceDeRham operator canonically as composition \( (\sharp\circ\laplaceDeRham\circ\flat) \).
Finally, we obtain
\begin{align}
  \laplaceDeRham \pb &= -\left(\laplaceRotRot  +  \laplaceGradDiv \right) \pb \formPeriod
\end{align}
for vector-valued functions (and 1-forms) $\pb$ (and $\alphab$).

\section{Numerical methods}\label{sec:NumericalMethods}
The growing interest in partial differential equations on surfaces is driven by various applications, but also by challenging numerical problems, which result from
the nonlinearity due to the underlying curved space. Various numerical methods have been developed to deal with these problems for scalar-valued surface partial differential equations. Finite element spaces are constructed on triangulated surface \cite{Dziuk1988,Dziuk2007a,Dziuk2007b}. These surface finite elements essentially allows to use the same concepts and tools as
in flat space \cite{AMDiS2007,Dziuk2013} and also the computational cost is comparable. The same holds for finite volume methods on quadrilateral grids on surfaces \cite{Calhoun2008}. Other
approaches consider an implicit representation of the surface, either through a level set description \cite{Bertalmio2001,Greer2006,Stoecker2008a,Dziuk2008}, within a diffuse interface
approximation \cite{Raetz2006,Raetz2007} or a closest point method \cite{Ruuth2008,Macdonald2008}. All these methods only require minimal information on the surface. All geometric
information is constructed solely through knowledge of the vertices of the discretization, or through the implicit description of a level set, phase field function, or point cloud. This has been proven to be sufficient and leads to efficient numerical methods also for complex physical problems \cite{Eilks2008,Lowengrub2009,Aland2011,Raetz2012,Aland2012,Nitschke2012,Stoop2015}.

For vector-valued surface partial differential equations the coupling between the equation and the geometry is much stronger and numerical methods which reduce the geometric information to a minimum might no longer be the most efficient. The literature on numerical methods for such problems is rare and mainly restricted to special surfaces, like the sphere. Here, spectral methods based on spherical harmonics expansions
are a popular tool \cite{Backus1966,Barrera1985,Freeden1994,Kostelec2000,Fengler2005,Freeden2009}. Another method which makes use of detailed geometric properties is an exterior calculus
approach \cite{Hirani2003,Desbrun2005,Arnold2006,Arnold2010}, which has recently also been applied to vector-valued surface partial differential equations, \eg, surface Navier-Stokes equations \cite{Mohamed2016,Nitschke2016}.

We will consider four different methods to solve the weak surface Frank-Oseen problem \eqref{eq:VectorPDE}, and
\eqref{eq:FormPDE}.
The first method is a Discrete Exterior Calculus (\DEC) formulation of equation \eqref{eq:FormPDE}, to be discussed in \autoref{sec:Dec}. Handling the penalty term requires an implementation of a pair of discrete equations for the dual vector and its hodge-dual variant and leads to a coupled system of primal-dual equations, which to the best of our knowledge has not been considered before in this context.
In \autoref{sec:SpectralMethod} the second method based on spherical harmonics (\SPH) is introduced. This approach expands $\pb$ in a spherical function basis, given as eigenfunctions of the \LaplaceDeRham operator. This results in a discrete set of equations for the expansion coefficients.
The third approach is the surface finite element method (\SFEM), to be explained in \autoref{sec:ParametricFiniteElements}. It relaxes the requirement of $\pb$ to be a tangential field, by introducing an additional penalty term that weakly enforces tangentiality. The vector field is represented in an Euclidean basis, leading to a system of scalar-valued surface PDEs. The representation of the \LaplaceDeRham operator in an Euclidean basis restricted to the tangent-plane by penalty terms is a new ansatz to discretize vector-valued surface PDEs.
The fourth method is the diffuse interface method (\DI), see \autoref{sec:DiffuseSurfaceApproximation}. It extends the domain to the embedding space $\R^3$, enforces tangentiality weakly and additionally restricts the differential operators to the surface using an approximation of a surface delta function. This leads to a system of coupled scalar-valued PDEs in a three dimensional domain and extends the established concept to vector-valued surface PDEs.

In the following section the time-discretization for the evolution problem is introduced. It is shared by all considered methods.

\subsection{\label{sec:time-discretization}Discretization in time}
Let $0 < t_0 < t_1 < \ldots$ be a sequence of discrete times with time step width $\tau_k := t_{k+1}-t_k$ in the $k$-th iteration. The fields $\pb^k(\xb)\equiv \pb(\xb,t_k)$
and $\alphab^k(\xb)\equiv \alphab(\xb,t_k)$, respectively, correspond to the time-discrete functions at $t_k$. Applying a semi-implicit Euler discretization to \eqref{eq:VectorPDE} and \eqref{eq:FormPDE} results in time discrete systems of equations as follows:
Let $\pb^0\in C(\surf;\,\Tangent\surf)$ be a given initial director field. For $k=0,1,2,\ldots$ find $\pb^{k+1}\in C^2(\surf;\,\Tangent\surf)$ s.t.
\begin{equation}\label{eq:time-discrete-system-p}
\frac{1}{\tau_k}\pb^{k+1} + K(\laplaceDeRham\pb^{k+1} + \shapeOperator^2\pb^{k+1}) + \Kn f(\pb^k,\pb^{k+1}) = \frac{1}{\tau_k}\pb^k \quad\text{in }\surf\formComma
\end{equation}
with $f(\pb^k,\pb^{k+1})$ a linearization of the non-linear term. In the methods \DEC, \SFEM, and \DI\ we consider a linear Taylor expansion  around $\pb^k$, see \eqref{eq:linearization_taylor}, and in the method \SPH\ we implement an explicit evaluation at the old time step $t_k$, see \eqref{eq:linearization_expl}:
\begin{align}
f^\text{Taylor}(\pb^k,\pb^{k+1}) &:= (\|\pb^k\|^2-1)\pb^{k+1} + 2\langle\pb^{k+1},\pb^k\rangle\pb^k - 2\|\pb^k\|^2\pb^k \label{eq:linearization_taylor}\\
f^\text{expl}(\pb^k,\pb^{k+1}) &:= \|\pb^k\|^2\pb^k - \pb^{k+1}\formPeriod \label{eq:linearization_expl}
\end{align}

The corresponding time discretization of the dual vector formulation \eqref{eq:FormPDE} is similar to \eqref{eq:time-discrete-system-p} utilizing the correspondence between vectors and dual vectors by the musical isomorphism $\flat$ for the initial condition: Let $\alphab_0:=\pb_0^\flat$ be given. For $k=0,1,2,\ldots$ find $\alphab_{k+1}\in\FormSpace(\surf)$ s.t.
\begin{equation}\label{eq:time-discrete-system-alpha}
\frac{1}{\tau_k}\alphab^{k+1} + K(\laplaceDeRham\alphab^{k+1} + \shapeOperator^2\alphab^{k+1}) + \Kn f(\alphab^k,\alphab^{k+1}) = \frac{1}{\tau_k}\alphab^k \quad\text{in }\surf\formPeriod
\end{equation}

\subsection{DEC}\label{sec:Dec}

For a Discrete Exterior Calculus the surface discretization is a simplicial complex \( \SC=\Vs\sqcup\Es\sqcup\Fs \)
containing sets of vertices \( \Vs \), edges \( \Es \), and (triangular) faces \( \Fs \).
The quantities of interest in our DEC discretization are 1-forms \( \alphab\in\FormSpace(\surf)=\Tangent^{*}\surf \).
We do not approximate the coordinate function of \( \alphab \) on a discrete set of points
or vertices, but rather introduce a finite set of degrees of freedom (DOFs) as integral values on the edges \( e\in\Es \),
\begin{align}\label{eq:discreteOneForm}
  \alpha_{h}(e) := \int_{\pi(e)}\alphab\formComma
\end{align}
with the gluing map \( \pi:\Es\rightarrow\surf \), which projects geometrically the edge \( e \) to the surface \( \surf \).
The mapping \(\alpha_{h}\in\FormSpace_{h}(\SC) \) is called the discrete 1-form of \( \alphab \), since \(  \alpha_{h}(e) \) approximates \( \alphab(e)\equiv\alphab(\eb) = \scalarprod{\pb,\eb} \) on an intermediate point \( \xi\in\pi(e)\subset\surf \),
where the edge vector \( \eb \) exists in \( \Tangent_{\xi}\surf|_{\pi(e)} \) by the mean value theorem.
Therefore, we approximate 1-forms on the restricted dual tangential space \( \Tangent_{\xi}\surf|_{\pi(e)} \),
which is a one dimensional vector space in \( \xi\in\surf \)
likewise the space of discrete 1-forms \( \FormSpace_{h}\left( \SC \right)|_{e} = \FormSpace_{h}\left( \{e\} \right)\)
restricted to the edge \( e \).
Hence, a discrete 1-form problem on surfaces leads to a one dimensional problem, like a scalar-valued problem.

The simplicial complex \( \SC \) is manifold-like, orientable and well-centered. For a detailed discussion of these requirements and general introduction to DEC,
see \cite{Hirani2003,Desbrun2005}.

Discrete linear differential operators composed of the exterior derivative \( \exd \) and the Hodge operator \( * \),
like the \LaplaceDeRham operator $\laplaceDeRham$, see \autoref{sec:Notation}, can be implemented by successively utilizing a discrete version of the Hodge operator
and the Stokes theorem for the exterior derivative, see \cite{Hirani2003}.
This procedure leads to a \DEC\ discretized Rot-Rot-Laplace $\laplaceRotRot_{h}$ and Grad-Div-Laplace $\laplaceGradDiv_{h}$.
For discrete 1-forms \( \alpha_{h}\in\FormSpace_{h}(\SC) \), sign mappings \( s_{\circ,\circ}\in\left\{ -1,+1 \right\} \),
volumes \( \left| \cdot \right| \), Voronoi cells \( \star v \), Voronoi edges \( \star e \), and the ``belongs-to'' relations \( \succ \) and \(
\prec\) we obtain
\begin{align}
  \laplaceRotRot_{h}\alpha_{h}(e) &= -\frac{\left| e \right|}{\left| \star e \right|} \sum_{\face\succ e} \frac{s_{\face,e}}{\left| \face \right|}
                                \sum_{\tilde{e}\prec\face} s_{\face,\tilde{e}}\, \alpha_{h}(\tilde{e}) \formComma\\
  \laplaceGradDiv_{h}\alpha_{h}(e) &= -\sum_{v\prec e} \frac{s_{v,e}}{\left| \star v \right|} \sum_{\tilde{e}\succ v}
                                s_{v,\tilde{e}} \frac{\left| \star\tilde{e} \right|}{\left| \tilde{e} \right|}\alpha_{h}(\tilde{e})
                                \formPeriod
\end{align}
Hence, in analogy to \autoref{sec:Notation}, we get the \DEC\ discretized \LaplaceDeRham operator
\( \laplaceDeRham_{h}\alpha_{h}(e) = -( \laplaceRotRot_{h}\alpha_{h}(e) + \laplaceGradDiv_{h}\alpha_{h}(e) ) \).
See \autoref{sec:DECLaplaceOperators} for details in notation and derivation of the \DEC\ operators.
The value for \( \laplaceDeRham_{h}\alpha_{h}(e) \) on an edge \( e \) is determined as a linear combination of few edge values
\(\alpha_{h}(\tilde{e}) \) in a \textit{proximate neighborhood} of \( e \),
\ie, it exists a vertex \( v \) that connects the edges \( e\succ v \) and \( \tilde{e}\succ v \).

Restricting the time-discrete evolution equation \eqref{eq:time-discrete-system-alpha} to the edges, using \eqref{eq:discreteOneForm}, leads to a system of equations for all edges $e\in\Es$:
\begin{equation}\label{eq:discrete_formPDE}
\frac{1}{\tau_k}\alpha^{k+1}_{h}(e)
                  + K\left(\laplaceDeRham_{h}\alpha^{k+1}_{h}(e) + \left(\shapeOperator^2\alphab^{k+1}\right)_h(e)\right) +
                  \Kn \left(f(\alphab^k,\alphab^{k+1})\right)_{h}(e) = \frac{1}{\tau_k}\alpha^k_{h}(e)\formComma
\end{equation}
with $(\shapeOperator^2\alphab^{k+1})_h(e) = \int_{\pi(e)}\shapeOperator^2\alphab^{k+1}$.
Using Taylor expansion \eqref{eq:linearization_taylor} in its covariant form, we obtain in the \( (k+1) \)-th time step for the non-linear term
\begin{align}\label{eq:taylor_endoform}
   \begin{aligned}
   \left(f(\alphab^k,\alphab^{k+1})\right)_{h}(e)
      &=  \int_{\pi(e)}(\|\alphab^k\|^2-1)\alphab^{k+1} + 2\langle\alphab^{k+1},\alphab^k\rangle\alphab^k - 2\|\alphab^k\|^2\alphab^k\\
      &= \int_{\pi(e)} \left( (\|\alphab^k\|^2-1)\idCoVec + 2\alphab^k\otimes\left( \alphab^k \right)^{\sharp} \right)\alphab^{k+1}\\
      &\phantom{=} -\int_{\pi(e)} 2\|\alphab^k\|^2\idCoVec\alphab^k \\
      &=: \left( \EndoTaylorLHS^{k}\alphab^{k+1} \right)_{h}(e) - \left( \EndoTaylorRHS^{k}\alphab^{k} \right)_{h}(e)\formComma
      \end{aligned}
\end{align}
with the identity map \( \idCoVec: \Tangent^{*}\surf \rightarrow \Tangent^{*}\surf \).

In the remaining section we discuss how to implement the norm $\|\alpha_{h}(e)\|$, the upcoming inner product $\langle \alpha_{h}(e), \alpha'_{h}(e) \rangle$ in the evaluation of the non-linear term,
and the endomorphisms  $\shapeOperator^2,\,\EndoTaylorLHS^{k}$, and $\EndoTaylorRHS^{k}$.

For the edge \( e_{0}:= e \) we choose another edge \( e_{1} \) in the proximate neighborhood of \( e_{0} \).
These two edges
define a vector space \( V_{\face} := \Span{\eb_{0},\eb_{1}} \) for the face \( \face\succ e_{0},e_{1} \)
at the contact vertex \( v\prec e_{0},e_{1} \).
A barycentric parametrization of \( V_{\face} \), regarding the basis vectors, results in a flat discrete metric
\begin{align}\label{eq:discreteMetric}
  \tensor{g} &= (\eb_{i}\cdot\eb_{j})de^{i}de^{j} \formComma
\end{align}
with the ordinary \( \R^{3} \) dot product and the canonical dual basis \( \left\{ de^{0}, de^{1} \right\} \),
which spans the flat vector space for covariant vectors.
So, we can construct a 1-form \( \alphav(e)\in \Tangent^{*}V_{\face}, \) which is constant on \( \face, \) by
\( \alphav(e) = \alpha_h(e_{i})de^{i} \).
Hence, if \( g^{ij} \) are the components of the inverse of the metric \eqref{eq:discreteMetric},
the square of the norm is given by
\begin{align}\label{eq:NormAddEdge}
  \|\alpha_h(e)\|^2 \equiv \left\| \alphav(e) \right\|^{2} = \alpha_h(e_{i})g^{ij}\alpha_h(e_{j}) \formPeriod
\end{align}
This norm strongly depends on the choice of the additional edge \( e_{1} \). Considering the Voronoi edge \( \star e \) (see \autoref{sec:DECNotations}),
which is not an edge in a pure simplicial sense, but a chain containing two edges orthogonal to $e$, one on the left face \( \face_{2}\succ e \)
and one on the right face \( \face_{1} \succ e \),
\ie, \( \star e = \star e|_{T_{1}} + \star e|_{T_{2}} \), leads to a stable pair of edges.
With a piecewise linear barycentric parametrization \( \boldsymbol{\gamma}:[0,1]\rightarrow \star e \) of the polygonal chain \( \star e \), with piecewise constant derivative $\|\boldsymbol{\gamma}'\|=\left| \star e \right|$ we can define the Voronoi edge vector
\begin{align*}
  \eb_{\star} &:= \boldsymbol{\gamma}' \in \Tangent\face_{1}\sqcup\Tangent\face_{2}\formPeriod
\end{align*}
This leads to the discrete metric in terms of the orthogonal basis \( \left\{ \eb, \eb_{\star} \right\} \) and
the dual basis \( \left\{ de, de^{\star} \right\} \),
\begin{align}\label{eq:PDMetric}
  \tensor{g} &= \left| e \right|^{2}\left( de \right)^{2} + \left| \star e \right|^{2}\left( de^{\star} \right)^{2}
\end{align}
and, with \( \alpha_h(\star e) \approx -\frac{\left| \star e \right|}{\left| e \right|}(*\alpha_h)(e) \) (see \cite{Hirani2003}),
the discrete (covariant) vector-valued 1-form
\begin{align}
  \alphav(e) &= \alpha_h(e)de + \alpha_h(\star e)de^{\star}
              \approx \alpha_h(e)de -\frac{\left| \star e \right|}{\left| e \right|}(*\alpha_h)(e)de^{\star}\formPeriod
\end{align}
The resulting vector spaces on all edges \( e\in\Es \) can be summarized as disjoint unions to
\( \PDT\Es := \bigsqcup_{e\in\Es}\bigsqcup_{\face\succ e} \Span{\eb, \eb_{\star}|_{\face}} \).
We call \( \big( \alpha_h, *\alpha_h \big):= \alphav: \Es \rightarrow \PDT^{*}\Es \)
a \textit{discrete primal-dual-1-form (PD-1-form)} with components \( \alpha_h \) and \( *\alpha_h \)
in \( \FormSpace_{h}(\SC) \).
Let \( \FormSpace_{h}(\SC;\PDT^{*}\Es) \) be the \textit{space of discrete PD-1-forms}.
All discrete PD-1-forms are uniquely defined
and depend only on the edge \( e \) and geometrical informations about it and its Voronoi edge.
Henceforward, we omit the argument \( e \) for a better readability.
The norm of \( \alphav=\alphav(e) \) is computed on all edges \( e \) with the discrete metric \eqref{eq:PDMetric} by
\begin{align}\label{eq:PDNorm}
  \left\| \alphav \right\|^{2} &= \frac{1}{\left| e \right|^{2}}\left( \alpha_h^{2} + (*\alpha_h)^{2} \right)
\end{align}
and the discrete inner product with another discrete PD-1-form \( \betav = \big(\beta_h, *\beta_h\big) \) is computed by
\begin{align}\label{eq:PDInner}
  \left\langle \alphav , \betav \right\rangle
      &= \frac{1}{\left| e \right|^{2}}\left( \alpha_h\beta_h + (*\alpha_h)(*\beta_h) \right) \formPeriod
\end{align}

The Hodge operator $*$ applied to \eqref{eq:time-discrete-system-alpha} results in the Hodge dual equation
\begin{align}\label{eq:DualFormPDE}
  \frac{1}{\tau_k} (*\alphab)^{k+1} + \K\left(\laplaceDeRham (*\alphab)^{k+1} + *\shapeOperator^{2}\alphab^{k+1}\right)
  + \Kn *f(\alphab^{k},\alphab^{k+1}) &= \frac{1}{\tau_k}(*\alphab)^k \formComma
\end{align}
where the identity \( *\laplaceDeRham = \laplaceDeRham* \) for the \LaplaceDeRham operator is used.
Restricting \eqref{eq:DualFormPDE} to the edges \( e\in\Es \), utilizing \eqref{eq:discreteOneForm}, and combining the result with \eqref{eq:discrete_formPDE}
leads to
\begin{multline}
  \frac{1}{\tau_k}\alphav^{k+1} + \K\left( \laplaceDeRham_h\alphav^{k+1} +
                \begin{bmatrix}
                  \left( \shapeOperator^{2}\alphab^{k+1} \right)_{h} \\
                  \left( *\shapeOperator^{2}\alphab^{k+1} \right)_{h}
                \end{bmatrix}\right)
            + \Kn
                  \begin{bmatrix}
                    \left(\EndoTaylorLHS^{k}\alphab^{k+1}\right)_{h} \\
                    \left(*\EndoTaylorLHS^{k}\alphab^{k+1}\right)_{h}
                  \end{bmatrix}\\
   = \frac{1}{\tau_k}\alphav^{k} + \Kn
                  \begin{bmatrix}
                    \left(\EndoTaylorRHS^{k}\alphab^{k}\right)_{h} \\
                    \left(*\EndoTaylorRHS^{k}\alphab^{k}\right)_{h}
                  \end{bmatrix} \text{\ in }\Es\formPeriod
\end{multline}

In \autoref{sec:DECLinOp} it is shown, how to approximate endomorphisms \( \tensor{M}:\Tangent^{*}\surf \rightarrow \Tangent^{*}\surf \) in a DEC-PD context,
so that
\begin{align}\label{eq:endomorphism_aprox}
  \begin{bmatrix}
                  \left( \tensor{M}\alphab \right)_{h} \\
                  \left( *\tensor{M}\alphab \right)_{h}
  \end{bmatrix}
      \approx \MPDt\cdot\alphav \text{\ in }\Es
\end{align}
with the mixed co- and contravariant \textit{discrete PD-(1,1)-Tensor} \( \MPDt \).
Evaluating the \( \R^{3} \) representation of the shape operator at
the midpoint of the edge \( e\in\Es \) projected to the surface, \ie, \( \shapeOperator^{2}(e):=\shapeOperator^{2}|_{\pi(c(e))}\in\R^{3\times 3} \),
utilizing \eqref{eq:EndoApprox}, results in a matrix form of the shape operator, applicable in \eqref{eq:endomorphism_aprox},
\begin{align}
   \shapeOperatorPDt^{2}(e) =
      \begin{bmatrix}
        \frac{\eb\cdot\shapeOperator^{2}(e)\cdot\eb}{|e|^{2}}
      &-\frac{\eb\cdot\shapeOperator^{2}(e)\cdot\eb_{\star}}{|e||\star e|} \\
       -\frac{\eb_{\star}\cdot\shapeOperator^{2}(e)\cdot\eb}{|e||\star e|}
      &\frac{\eb_{\star}\cdot\shapeOperator^{2}(e)\cdot\eb_{\star}}{|\star e|^{2}}
      \end{bmatrix}\formPeriod
\end{align}
Similarly, with \eqref{eq:PDNorm} and \( \idCoVec^{\flat} = \tensor{g} \), considering the discrete metric, we get
\begin{align}
  \EndoTaylorRHSPDt^{k}(e) &= 2\left\| \alphav^{k}(e) \right\|^{2} \begin{bmatrix} 1 & 0 \\ 0 & 1\end{bmatrix} \formComma\\
  \EndoTaylorLHSPDt^{k}(e) &= \left( \left\| \alphav^{k}(e) \right\|^{2} - 1 \right) \begin{bmatrix} 1 & 0 \\ 0 & 1\end{bmatrix}
                        +2\begin{bmatrix}
                          \frac{\alpha^{k}_{h}(e)\alpha^{k}_{h}(e)}{|e|^{2}} & - \frac{\alpha^{k}_{h}(e)\alpha^{k}_{h}(\star e)}{|e||\star e|} \\
        - \frac{\alpha^{k}_{h}(e)*\alpha^{k}_{h}(\star e)}{|e||\star e|} & \frac{\alpha^{k}_{h}(\star e)\alpha^{k}_{h}(\star e)}{|\star e|^{2}}
                        \end{bmatrix} \\
         &\approx \left( \left\| \alphav^{k}(e) \right\|^{2} - 1 \right) \begin{bmatrix} 1 & 0 \\ 0 & 1\end{bmatrix}
                +\frac{2}{|e|^{2}}\begin{bmatrix}
                      \alpha^{k}_{h}(e)\alpha^{k}_{h}(e)  & \alpha^{k}_{h}(e)\left( *\alpha \right)^{k}_{h}(e) \\
                      \alpha^{k}_{h}(e)\left( *\alpha \right)^{k}_{h}(e)  & \left( *\alpha \right)^{k}_{h}(e)\left( *\alpha \right)^{k}_{h}(e)
            \end{bmatrix}\notag\\
         &=: \tilde{\EndoTaylorLHSPDt}^{k}(e)\formPeriod\notag
\end{align}
Finally, with the discrete inner product \eqref{eq:PDInner},
\( \EndoTaylorRHSPDt^{k}\cdot\alphav^{k} = 2\left\| \alphav^{k} \right\|^{2}\alphav^{k} \),
and
\begin{align}
  \tilde{\EndoTaylorLHSPDt}^{k}\cdot\alphav^{k+1} &= \left( \left\| \alphav^{k} \right\|^{2} - 1 \right)\alphav^{k+1} + 2\left\langle \alphav^{k+1}, \alphav^{k} \right\rangle\alphav^{k}
    \text{\ in }\Es\formComma
\end{align}
the introduced Taylor linearization of $f$, \ie, $f^\text{Taylor}(\alphav^{k}, \alphav^{k+1})$, is found.

This results in a series of time-discrete linear DEC-PD problems: For $k=0,1,2,\ldots$, and a given initial value \( \alphav^{0} \), find \( \alphav^{k+1}\in\FormSpace_{h}(\SC;\PDT^{*}\Es) \) s.t.
\begin{align}\label{eq:DECPDLin}
  \frac{1}{\tau_{k}}\alphav^{k+1} + \K\left(\laplaceDeRham_h\alphav^{k+1} + \shapeOperatorPDt^{2}\cdot\alphav^{k+1} \right)
      + \Kn f^\text{Taylor}(\alphav^{k}, \alphav^{k+1}) &= \frac{1}{\tau_{k}}\alphav^{k} &\text{ in }\Es \formPeriod
\end{align}

These stationary problems can be implemented\footnote{For a software framework, see also the discretization library Dune-DEC \cite{DuneDEC}.} by assembling a matrix and vector for the components $\alpha_h(e)$ and $(*\alpha_h)(e)$ on edges $e\in\Es$.
The resulting linear system is solved with the TFQMR method, see \cite{Freund1993}.

Many conceivable ways exist to interpolate the initial condition
 \( \alphav^{0}\in\FormSpace_{h}(\SC;\PDT^{*}\Es) \), with \( \alphav^{0}= [\alpha^{0},*\alpha^{0}] \),
 from a given vector field \( \pb^{0}\in T\surf \).
We assume that the simplicial complex and its polytope \( |\SC| \) are immersed in
a sufficiently small neighborhood $\Omega_\delta$ of the surface, so that the
initial condition $\pb^{0}$ can be smoothly extended.

Given such an extension \( \widetilde{\pb}^{0} \) of an initial vector field \( \pb^{0}  \) we
 can choose the intersection point \( c(e) \) of an edge \( e\in\Es \) and \( \star e \) for approximating the integral
 expressions, \ie, let the edge \( e \) be given so that it points from the vertex \( v_{1} \) to the vertex \( v_{2} \)
 and the dual edge \( \star e \) from the circumcenter \( c(\face_{1}) \) to \( c(\face_{2}) \), then we obtain
 \begin{align}
   \alpha^{0}(e) &= \int_{\pi(e)}\left( \pb^{0} \right)^{\flat}
                 \approx \int_{0}^{1}\widetilde{\pb}^{0}({c(e)})\cdot\eb\,\dif\lambda
                 = \widetilde{\pb}^{0}({c(e)})\cdot(v_{2} - v_{1})\\
   (*\alpha)^{0}(e) &\approx -\frac{\left| e \right|}{\left| \star e \right|} \int_{\star\pi(e)}\left( \pb^{0} \right)^{\flat}
                 \approx -\frac{\left| e \right|}{\left| \star e \right|} \int_{0}^{1}\widetilde{\pb}^{0}({c(e)})\cdot\eb_{\star}\,\dif\lambda\\
                &= -\frac{\left| e \right|}{\left| \star e \right|} \widetilde{\pb}^{0}({c(e)})\cdot(c(\face_{2}) - c(\face_{1}))
                \formPeriod \notag
 \end{align}

On the other hand, if \( \pb^{0} \) arise from the gradient of a scalar function \( f:\surf\to\R \),
 \ie, \( \pb^{0} = \Grad f = (\exd f)^{\sharp}\),
 we obtain for a smooth extension \( \widetilde{f} \) of $f$:
 \begin{align}
   \alpha^{0}(e) &= \int_{\pi(e)} \exd f = f(v_{2}) - f(v_{1}) \\
   (*\alpha)^{0}(e) &\approx -\frac{\left| e \right|}{\left| \star e \right|} \int_{\star\pi(e)} \exd f \\
                 &\approx -\frac{\left| e \right|}{\left| \star e \right|}
                             \left( \widetilde{f}(c(\face_{2})) -  \widetilde{f}(c(\face_{1}))\right) \formComma \notag
 \end{align}
  utilizing Stoke's theorem.
 
\subsection{Spectral method}\label{sec:SpectralMethod}
In this section we restrict our consideration to spherical surfaces $\surf=\Sp$ parametrized
by $\lu\in[0,\pi]$ and $\lv\in[0,2\pi)$, \ie, the co-latitude and azimuthal coordinates, respectively.
So each point $\mathbf{x}_{\Sp}\in\Sp$ can be written as $\mathbf{x}_{\Sp}(\lu,\lv)=\sin(\lu)\cos(\lv)\EuBase{x} + \sin(\lu)\sin(\lv)\EuBase{y} + \cos(\lu)\EuBase{z}$.
Based on the observation that the tangential part of a spherical vector field
can be split into a curl-free and a divergence-free field by using derivatives of
scalar fields, an efficient numerical methods can be constructed. The Helmholtz
decomposition theorem \cite{Freeden2009} states that every continuously differentiable spherical tangent
vector field $\vect{f}:\surf\rightarrow \Tangent\surf$ can be represented by uniquely
determined scalar functions $f_1,f_2\in\Csurf{1}$ as
\[
\vect{f}(\xb) = \Grad f_1(\xb) + \Rot f_2(\xb).
\]
An efficient solution method for linear
surface PDEs on the sphere is based on a spectral expansion of
the objective scalar functions $f\in\LS$ in the spherical harmonics $Y_l^m:\surf\rightarrow\mathbb{C}$,
$(l,m)\in\mathcal{I}_\infty$ with
$\mathcal{I}_N := \{ (l,m)\,:\,0\leq l\leq N, -l\leq m\leq l\}$,
which build an $\LS$-orthonormal system of
eigenfunctions of the Laplace-Beltrami operator $\laplaceBeltrami$, \ie,
\begin{equation}
\laplaceBeltrami Y_l^m = \Delta_{lm} Y_l^m\;\text{ with }\Delta_{lm}:=-l(l+1),\quad\text{ for }(l,m)\in\mathcal{I}_\infty
\end{equation}
and $\LtwoProd{Y_l^m,\, Y_{l'}^{m'}} = \delta_{ll'}\delta_{mm'}$, cf. \cite{Hesthaven2007, Backofen2011}.
Due to the symmetries of the sphere, analytic representations of $Y_l^m$ can be found in terms of Associated Legendre polynomials.
This allows for an efficient evaluation of the basis functions.

A scalar function $f\in L^2(\surf)$ can be represented in the series expansion
\begin{equation}\label{eq:expansion_f}
f(\lu, \lv) = \sum_{l=0}^\infty\sum_{m=-l}^l f_{lm} Y_l^m(\lu, \lv)
\end{equation}
with expansion coefficient $f_{lm} = \LtwoProd{f,\,Y_l^m}$.

Taking the gradient
and curl of the spherical harmonics, an expansion for tangential vector fields can
be constructed. Therefore, we introduce two vector spherical harmonics
$\boldsymbol{y}_{lm}^{(1)},\,\boldsymbol{y}_{lm}^{(2)}$ as
\begin{align}\label{eq:vector_spherical_harmonics}
\begin{aligned}
  \boldsymbol{y}_{lm}^{(1)}(\lu, \lv) &:= N_{lm}\Grad Y_l^m(\lu, \lv)\,,\\
  \boldsymbol{y}_{lm}^{(2)}(\lu, \lv) &:= N_{lm}\Rot Y_l^m(\lu, \lv)
\end{aligned}
\end{align}
with normalization constants $N_{lm} = (-\Delta_{lm}^{-1})^{1/2}$. These functions are normalized in
such a way, that they build again an $\Ltwo{\surf}{\Tangent\surf}$-orthonormal system of eigenfunctions
of a Laplace operator, namely the spherical \LaplaceDeRham operator \cite{Freeden1994,Freeden2009}, \ie,
\begin{equation}
\laplaceDeRahm\boldsymbol{y}^{(i)}_{lm} = -\Delta_{lm} \boldsymbol{y}^{(i)}_{lm},\quad\text{ for }i=1,2,\; (l,m)\in\mathcal{I}_\infty
\end{equation}
and $\left(\boldsymbol{y}^{(i)}_{lm},\, \boldsymbol{y}^{(j)}_{l'm'}\right)_{\Ltwo{\surf}{\Tangent\surf}} = \delta_{ij}\delta_{ll'}\delta_{mm'}$.

A series expansion of a tangent vector field $\vect{f}\in \Ltwo{\surf}{\Tangent\surf}$, based
on the expansion of scalar fields \eqref{eq:expansion_f} and the gradient and curl
basis representation \eqref{eq:vector_spherical_harmonics}, can thus be written as
\begin{equation}\label{eq:expansion_vec_f}
\vect{f}(\lu,\lv) = \sum_{i=1}^2\sum_{l=0}^\infty\sum_{m=-l}^l f^{(i)}_{lm} \boldsymbol{y}^{(i)}_{lm}(\lu, \lv)
\end{equation}
with expansion coefficients $f^{(i)}_{lm} = \left(\vect{f},\,\boldsymbol{y}^{(i)}_{lm}\right)_{\Ltwo{\surf}{\Tangent\surf}}$.
In the following we use the notation $\fv_{lm} := \big(f^{(1)}_{lm},\,f^{(2)}_{lm}\big)$ to
denote the pair of coefficients.

The spherical harmonics method is based
on the idea to approximate any scalar function
\eqref{eq:expansion_f} and vector-valued function \eqref{eq:expansion_vec_f} by truncated expansions with band-width $l_{\max}=: N$.
Therefore, we introduce the space of spherical vector polynomials
\[
\vec{\Pi}_N(\surf) := \Big\{ \mathbf{f} = \sum_{i=1}^2\sum_{l=0}^N\sum_{m=-l}^l f^{(i)}_{lm} \boldsymbol{y}^{(i)}_{lm} \Big\}.
\]
The evaluation of expansion coefficients, in other words, the calculation of the $L^2$ inner product, is implemented by approximating the
integral by an appropriate quadrature rule. Let $\Vs = \{\Xb_k = (\lu_k, \lv_k)\}$ be a set of quadrature points on the sphere and $\{w_k\}$ the corresponding
quadrature weights. We introduce the discrete $L^2$ inner product:
\[
\left(\vect{f},\,\boldsymbol{y}^{(i)}_{lm}\right)_{h, \Ltwo{\surf}{\Tangent\surf}} := \sum_k w_k \Scalarprod{\vect{f}(\lu_k,\lv_k),\,\bar{\boldsymbol{y}}^{(i)}_{lm}(\lu_k,\lv_k)}
\]

In order to derive an equation for the expansion coefficients $\pv_{lm}$ of $\pb\in\vec{\Pi}_N(\surf)$
in terms of a Galerkin approach, see, \eg, \cite{Hesthaven2007},
we require the residual $\boldsymbol{r}$ of the differential equation \eqref{eq:VectorPDE},
\[
\boldsymbol{r} := \frac{1}{\tau_k}(\pb^{k+1} - \pb^k) + K(\laplaceDeRham\pb^{k+1} + \shapeOperator^2\pb^{k+1}) + \Kn f(\pb^k,\pb^{k+1})\,,
\]
to be orthogonal to the basis of $\vec{\Pi}_N(\surf)$ w.r.t. the $L^2$ inner product, \ie,
\begin{equation}\label{eq:sph_residual_equation}
\left(\boldsymbol{r},\, \boldsymbol{y}^{(i)}_{lm}\right)_{\Ltwo{\surf}{\Tangent\surf}} = 0,\quad\text{ for }i=1,2,\; (l,m)\in\mathcal{I}_N\formPeriod
\end{equation}
The shape operator on $\Sp$ simplifies to the surface identity, \ie, $\shapeOperator=-\ProjectSurf$.
With $f:=f^\text{expl}$, this term can be evaluated in discrete grid points on the sphere
rather than by forming convolution sums of the coefficients, see \cite{Boyd2001}. Therefore, let
the non-linear term $\mathbf{f}^k:=\norm{\pb^k}^2\pb^k$ at time step $t_k$ be expanded in the space $\vec{\Pi}_N(\surf)$
with expansion coefficients $\underline{\mathbf{f}}_{lm}^k$.

By requiring the new time step solution $\pb^{k+1}$ to be an element of $\vec{\Pi}_N(\surf)$, we can insert the truncated expansion of the solution into the residual equation \eqref{eq:sph_residual_equation}. Utilizing the property that the \LaplaceDeRham operator is the eigen-operator of the basis functions results in an equation for the
expansion coefficients directly. Finally, the time step procedure for the spherical harmonics approach reads: Let
\[
\pv^{0,(i)}_{lm} = \left(\pb^0,\,\boldsymbol{y}^{(i)}_{lm}\right)_{h, \Ltwo{\surf}{\Tangent\surf}}\quad\text{ for }i=1,2,\;(l,m)\in\mathcal{I}_N
\]
be the expansion coefficients for the initial solution.
For $k=0,1,2,\ldots$ \begin{enumerate}
\item Evaluate $\mathbf{f}^k(\Xb):=\norm{\pb^k(\Xb)}^2\pb^k(\Xb)$ for all $\Xb\in\Vs$.
\item Calculate $f^{k,(i)}_{lm} = \left(\vect{f}^k,\,\boldsymbol{y}^{(i)}_{lm}\right)_{h, \Ltwo{\surf}{\Tangent\surf}}$
        for \(i=1,2,\;(l,m)\in\mathcal{I}_N\).
\item Solve
\[
\frac{1}{\tau_k}\pv^{k+1}_{lm} - K\Delta_{lm}\pv^{k+1}_{lm} + (K - \Kn)\pv^{k+1}_{lm} = \frac{1}{\tau_k}\pv^{k}_{lm} - \Kn\fv^k_{lm},\quad\forall (l,m)\in \mathcal{I}_N
\]
to be understood component-wise.
\item Evaluate \eqref{eq:expansion_vec_f} with coefficients $\pv^{k+1}_{lm}$ to get $\pb^{k+1}$.
\end{enumerate}

The discrete spherical harmonics transform, that is, the evaluation of \eqref{eq:expansion_f} for a band-width $N$, can be split up into a discrete Fourier transform, realizable by a fast Fourier transform, and discrete Legendre transforms, implemented thanks to discrete cosine transforms \cite{Kunis2003} or a fast multipole method \cite{Suda2002}. The inverse transform, \ie, the calculation of the expansion coefficients, may be realized by the Gauss-Legendre algorithm. There, the integral is replaced by a Gauss-Legendre quadrature rule with Gauss nodes and weights in latitudinal direction \cite{Schaeffer2013}. Therefore, the spherical coordinate space is discretized by the set of vertices
\[
\Vs:=\left\{\Xb(\lu_i,\lv_j)\,:\,0\leq i< N_{\lu}\,,\,0\leq j < N_{\lv}\right\}\formComma
\]
with $\lu_i$ Gauss nodes in $[0,\pi]$ and $\lv_j$ equally distributed in $[0,2\pi)$. To respect the sampling theorem, we have chosen $N_{\lu}>N$ and $N_{\lv}>2 N$. Therewith, the coefficients of the non-linear term are only approximated, since $\textbf{f}^k$ is not in $\vec{\Pi}_N(\surf)$ for $\pb^k\in\vec{\Pi}_N(\surf)$.

Finally, the discrete vector harmonic transform can be implemented by two scalar transforms, see, \eg, \cite{Kostelec2000}. Thus, the complexity of the transform is dominated by the scalar transform that can be realized in $\mathcal{O}(N^2\log{N})$ \cite{Suda2002}.
 
\subsection{Surface finite elements}\label{sec:ParametricFiniteElements}
We consider a reformulation of $\EE$ and dynamic equation \eqref{eq:VectorPDE} suitable for a component-wise surface finite element approximation. To do so, we extend $\EE$ to a domain of vector-valued functions $\pExt:\surf \to \Tangent\R^3$ and penalize any energy contributions by normal components $\pExt \cdot \surfNormal \neq 0$ with a penalty factor $\Kt \gg 1$. The previously introduced \LaplaceDeRham operator has been defined as a differential operator on sections of tangent bundles. This needs to be extended to $\R^3$ vector fields. In a first step we use the surface projection $\ProjectSurf$ introduced in $\eqref{eq:surfaceProjection}$ and a result from \cite{Duduchava2006} to express $\Div \pb$ by $\Div\pExt$, \ie,
\begin{align}
\Div \pb = \Div \left(\ProjectSurf \pExt \right) = \underbrace{\nabla \cdot \pExt - \surfNormal \cdot (\nabla \pExt \cdot \surfNormal)}_{= \Div \pExt} - \meanCurvature\left( \pExt \cdot \surfNormal \right)\formComma
\end{align}
where $\meanCurvature=\Div\surfNormal$ denotes the mean curvature of $\surf$. Note that the curl of a vector field reduces to the curl of its tangential part, \ie, $\Rot \pb = \Rot \pExt$.  Further, we apply a decomposition of $\pExt = \pb + \surfNormal \left( \pExt \cdot \surfNormal \right)$ and $\qExt = \qb + \surfNormal \left( \qExt \cdot \surfNormal \right)$ to express the $L^2$ inner product of $\laplaceDeRham \pb$ and $\qb$ in terms of $\pExt$ and $\qExt$ (for details see \autoref{sec:IntegralTheorems}),
\begin{align*}
 \int_{\surf} \langle \laplaceDeRham \pb, \qb \rangle \dS = & \int_{\surf} \left(\Div \pExt \right) \left(\Div \qExt \right) + \left(\Rot \pExt \right)  \left( \Rot\qExt \right)\dS + \int_{\surf} \meanCurvature^2 \left( \pExt \cdot \surfNormal \right) \left( \qExt \cdot \surfNormal \right) \dS \\
  - & \int_{\surf} \meanCurvature\left( \left( \qExt \cdot \surfNormal \right)\left(\Div \pExt \right) + \left( \pExt \cdot \surfNormal \right) \left(\Div\qExt \right)  \right) \dS \nonumber \formPeriod
\end{align*}
In order to neglect the terms involving normal components $\left( \pExt \cdot \surfNormal \right)$ and $\left( \qExt \cdot \surfNormal \right)$, the penalty term $\frac{\Kt}{2} \left(\pExt \cdot \surfNormal \right)^2$ is added to the energy $\EE$. The functional derivative of this contribution results in a symmetric term
\begin{align}
 \label{eqn:parFEMsymNormal}
 \int_{\surf} \frac{\Kt}{2} \Scalarprod{\Fdif{\left(\pExt \cdot \surfNormal \right)^2 }{\pExt}[\pExt],\qExt} \dS = \int_{\surf}\Kt \left(\pExt \cdot \surfNormal \right) \left( \qExt \cdot \surfNormal \right ) \dS
\end{align}
leading in the context of a minimization process to $\left( \pExt \cdot \surfNormal \right) \rightarrow 0$ and $\left( \qExt \cdot \surfNormal \right) \rightarrow 0$ as $\Kt \rightarrow \infty$. As a result, we obtain an approximation of the \LaplaceDeRham operator for finite $\Kt$ by
\begin{align}\label{eq:parFEMdefLaplaceTilde}
\begin{aligned}
  \int_{\surf} \langle \laplaceDeRham \pb, \qExt \rangle \dS \approx & \int_{\surf} \left(\Div \pExt \right) \left(\Div \qExt \right) + \left(\Rot \pExt\right) \left(\Rot \qExt \right)\dS \\
 & = \int_{\surf} \underbrace{-\left[ \Grad \left(\Div \pExt \right)  + \Rot \left( \Rot \pExt \right) \right]}_{=\laplaceDeRhamTilde \pExt} \cdot\, \qExt \dS \formPeriod
\end{aligned}
\end{align}
A brief numerical study justifying this approach can be found in \autoref{sec:convApproxLaplaceDeRham}.
With this established, we formulate the extended weak surface Frank-Oseen energy for $\pExt \in \HdrExt{\surf}{\R^3}$ as:
\begin{align}
\begin{aligned}
  \F{\Kn,\Kt}^\surf[\pExt]  & =  \int_{\surf}  \frac{\K}{2}\left[\left( \Div \pExt \right)^2 + \left( \Rot \pExt \right)^2  + \| \shapeOperator \cdot \pExt \|^{2} \right] \dS \\
& + \int_{\surf} \frac{\Kn}{4}\left( \| \pExt \|^2 - 1 \right)^2 + \frac{\Kt}{2} \left(\pExt \cdot \surfNormal \right)^2 \dS \formPeriod
\end{aligned}
\end{align}

A straightforward first variation of the energy leads to the associated equation
\begin{align}
 \dt \pExt + \K \left(\laplaceDeRhamTilde\pExt + \shapeOperator^2 \pExt \right) + \Kt \left(\surfNormal  \cdot \pExt \right) \surfNormal + \Kn \left( \| \pExt \|^2 - 1 \right) \pExt = 0 & \text{ in } \surf\times(0,\infty)
\end{align}
with the initial condition $\pExt(t=0) = \pb^{0} \in \Tangent \surf$. Using the vector space property of the extended variational space $\HdrExt{\surf}{\R^3}$ we split the vector-valued variational problem into a set of component-wise scalar variational problems\footnote{Here, we use lower indices to denote the components of a vector, not to mix up with the covariant indices used in the context of differential geometry.}. Therefore, let $\qExt$ be decomposed as
\begin{align}
 \HdrExt{\surf}{\R^3} \supseteq \left[\Hi{\surf}\right]^3 \ni \qExt = \sum_{i=1}^3 \extend{\q}_i \vect{e}_i, \quad \extend{\q}_i \in \Hi{\surf}\formComma
\end{align}
with $\{\vect{e}_i\}_i$ the Euclidean basis of $\R^3$. We obtain a set of coupled variational problems for $\extend{p}_i \in \Ltwo{0,\infty}{\Hi{\surf}}$
\begin{multline}
\label{eq:parFEMvarProblem}
 \int_{\surf} \dt \extend{\p}_{i} \extend{\q}  \dS
 +\int_{\surf}  \K \left[\left( \Div \pExt \right) \left(\Grad \extend{\q} \, \right)_i +  \left(\Rot \pExt \right) \left( \Rot(\extend{\q}\,\vect{e}_i) \right) + \left(\shapeOperator^2 \cdot \pExt\right)_i \extend{\q} \right]\dS\\
+ \int_{\surf}\Kt \left(\surfNormal  \cdot \pExt \right) \surfNormalI_{i} \extend{\q} + \Kn \left( \| \pExt \|^2 - 1 \right) \extend{\p}_{i} \extend{\q} \dS  = 0, \quad \forall \; \extend{\q} \in \Hi{\surf} \; \forall \; t \in \left( 0, \infty\right)
\end{multline}
for $i = 1, \ldots, 3$. To solve this set of variational problems, we have implemented the time-discretization introduced in \autoref{sec:time-discretization}. The tangential penalty term is evaluated at the new time step $t_{k+1}$ and the non-linear term is linearized using the expression $f^\text{Taylor}$. For the discretization in space, we apply the surface finite element method for scalar-valued PDEs \cite{Dziuk1988,Dziuk2007a,Dziuk2007b} for each component.
Therefore, the surface $\surf$ is discretized by a conforming triangulation $\triangulation$, given as the union of simplices in a simplicial complex, \ie,
\[
\triangulation := \bigcup_{\sigma\in\SC} \sigma\formPeriod
\]
We use globally continuous, piecewise linear Lagrange elements
\[
\mathbb{V}_h(\triangulation) = \left\{ v_h \in C^0(\triangulation) \,:\, v_h|_T \in \mathbb{P}^1, \, \forall \, T \in \Fs \right\}
\]
as trial and test space for all components $\extend{\p}_{i}$ of $\pExt$, with $\Fs$ the set of triangular faces.

The resulting discrete problem reads: For $k=0,1,2,\ldots$ find $\extend{\p}_{i}^{k+1}\in\mathbb{V}_h(\triangulation)$ s.t.
\begin{multline}
 \frac{1}{\tau_k} \int_{\surf_h} \extend{\p}_{i}^{k+1}\extend{\q}\dS + \K  \int_{\surf_h} \Div \pExt^{k+1} \left(\Grad \extend{\q}\, \right)_i + \Rot \pExt^{k+1} \Rot (\extend{\q} \vect{e}_i) + \left(\shapeOperator^2
 \cdot \pExt\right)_i \extend{\q} \dS \\
 + \Kt \int_{\surf_h} \surfNormal\cdot\pExt^{k+1} \surfNormalI_{i} \extend{\q} \dS + \Kn \int_{\surf_h} \left( \|\pExt^{k}\|^2 - 1 \right) \extend{\p}_{i}^{k+1}\extend{\q} + 2 \extend{\p}_{i}^{k} \pExt^{k}\cdot \pExt^{k+1} \extend{\q} \dS \\
 = \frac{1}{\tau_k} \int_{\surf_h} \extend{\p}_{i}^{k}\extend{\q}\dS + 2\Kn \int_{\surf_h} \|\pExt^{k}\|^2 \extend{\p}_{i}^{k} \extend{\q} \dS,  \quad \forall\,\extend{\q} \in \mathbb{V}_h(\triangulation) \;
\end{multline}
for $i = 1, \ldots, 3$. To assemble and solve the resulting system we use the FEM-toolbox AMDiS \cite{AMDiS2007,AMDiS2015} with domain decomposition on 8 processors. As linear solver we have used a restarted GMRES method with a restart cycle of 30, modified Gram-Schmidt orthogonalization, and a block Jacobi preconditioner with ILU($0$) local solver on each partition.

\subsection{Diffuse interface approximation}\label{sec:DiffuseSurfaceApproximation}
Based on the penalty formulation, described in \autoref{sec:ParametricFiniteElements}, we formulate a diffuse interface approximation following the general treatment introduced in \cite{Raetz2006}. We use a simple (e.g. box like) embedding domain $\surf \subset \domain \subset \mathbb{R}^3$ and describe the surface
as the $1/2$ levelset of a phase-field variable $\phi$ defined on $\domain$:
\begin{align}
 \phi(\xb) = \frac{1}{2}\left( 1 - \tanh \left( \frac{3}{\varepsilon} d_\surf(\xb)\right) \right)\,,
\end{align}
with interface thickness $\varepsilon$ and $d_\surf(\xb)$ a signed-distance function. This gives an approximation of
the surface delta function
\begin{align}
 \delta_{\surf} \simeq \frac{36}{\varepsilon} \phi^2(\phi-1)^2 = \doubleWell(\phi).
\end{align}
In this diffuse interface framework we consider vector fields $\extendDomain{\pb}:\Omega\to \Tangent\R^3$ extended from the surface to the embedding domain $\domain$. The outward pointing surface normals are extended smoothly to $\domain$ by using $\nExt = \nabla \phi / \|\nabla \phi\|$ and the shape operator in the embedding domain is defined in terms of this extended normal, \ie, $\extendDomain{\shapeOperator}_{ij} = -\left[\ProjectSurf \nabla \nExt_j\right]_i $.

Considering the diffuse interface approximation of the extended weak surface Frank-Oseen energy
\begin{align}
 \F{\Kn,\Kt}^\domain[\extendDomain{\pb}]
 & = \int_{\domain} \frac{\K}{2} \doubleWell(\phi) \left[ \left( \Div \extendDomain{\pb} \right)^2 + \left(\Rot \extendDomain{\pb}  \right)^2 + \|\extendDomain{\shapeOperator} \cdot \extendDomain{\pb}\|^2 \right] \dS \\
 & +  \int_{\domain} \frac{\Kn}{4} \doubleWell(\phi) \left( \| \extendDomain{\pb} \|^2 - 1 \right)^2 + \frac{\Kt}{2} \doubleWell(\phi) \left(\extendDomain{\pb} \cdot \nExt \right)^2 \dS \nonumber
\end{align}
with $\extendDomain{\pb} \in \HdrDiffuse{\domain}{\R^3}$, we obtain, by straightforward first variation, the \( L^{2} \)-gradient flow formulation
\begin{equation}\label{eqnDiffuseEqNonlin}
 \doubleWell(\phi)\dt \extendDomain{\pb} + \K \laplaceDeRhamDiffuse \extendDomain{\pb}
 +  \doubleWell(\phi)\left[\K \left(\extendDomain{\shapeOperator}^2 \cdot \extendDomain{\pb}\right) + \Kt  \left(\extendDomain{\surfNormal} \cdot \extendDomain{\pb} \right) \extendDomain{\surfNormal} + \Kn \left( \| \extendDomain{\pb} \|^2 - 1 \right) \extendDomain{\pb}\right] = 0
\end{equation}
 in $\domain\times(0,\infty)$. Here, we have introduced the diffuse interface Laplace-deRham operator $\laplaceDeRhamDiffuse$ by
 \begin{align}\label{eq:defDiffuseLaplaceDeRham}
  \laplaceDeRhamDiffuse \extendDomain{\pb} := -\left[ \nabla \left( \doubleWell(\phi) \nabla \cdot \extendDomain{\pb} \right) + \nExt \times \nabla \left( \doubleWell(\phi) \nabla \cdot \left( \extendDomain{\pb} \times \nExt \right) \right) \right]\formPeriod
 \end{align}
As initial condition we set $\extendDomain{\pb}(t=0) = \extendDomain{\pb}^{0}$ in $\domain$ such that $\extendDomain{\pb}^{0} |_{\surf} = \pb^0$. As boundary condition we specify
 \[
 \nabla \extendDomain{\p}_{i} \cdot \nb = 0,\quad\text{ on }\partial\domain\times(0,\infty)\formComma
 \]
for $i=1,\hdots,3$, where $\nb$ denotes the outward pointing normal of $\partial \domain$. For $\Omega$ big enough, the condition on the outer boundary does not influence the solution on the surface. Finally, we obtain a set of coupled variational problems for $\extendDomain{p}_i \in \Ltwo{0,\infty}{\Hi{\domain}}$
\begin{align}\label{eq:diffInterfaceVarProb}
 & \;\;\;\; \int_{\domain} \doubleWell(\phi) \dt \extendDomain{\p}_{i} \extendDomain{\q}  \dV \\
 &+\int_{\domain}  \K \,  \doubleWell(\phi)\left[  \left( \nabla \cdot \extendDomain{\pb} \right) \partial_i \extendDomain{\q} +   \nabla \cdot \left( \extendDomain{\pb} \times \nExt\right) \nabla \cdot \left( \extendDomain{\q} \vect{e}_i  \times \nExt\right) +  \left(\extendDomain{\shapeOperator}^2 \cdot \extendDomain{\pb}\right)_i \extendDomain{\q} \right]\dV \nonumber \\
& + \int_{\domain}\Kt \doubleWell(\phi) \left(\nExt  \cdot \extendDomain{\pb} \right) \extendDomain{\surfNormalI}_{i} \extendDomain{\q} + \Kn \doubleWell(\phi) \left( \| \extendDomain{\pb} \|^2 - 1 \right) \extendDomain{\p}_{i} \extendDomain{\q} \dV \nonumber \\
& = 0 \quad \forall \; \extendDomain{\q} \in \Hi{\domain} \; \forall \; t\in\left(0, \infty \right)\formComma\nonumber
\end{align}
for $i=1,\hdots,3$.

The definition of $\laplaceDeRhamDiffuse$ in \eqref{eq:defDiffuseLaplaceDeRham} is motivated by the component-wise formulation of $\laplaceDeRhamTilde$ in combination with the diffuse approximations of surface differential operators for scalar functions $f:\surf\to\R$ with smooth extension $\extendDomain{f}:\Omega\to\R$. In this framework we have the following convergence results:
\begin{align*}
\lim\limits_{\varepsilon \to 0}\int_{\domain }\doubleWell(\phi)\extendDomain{f}\,\extendDomain{\q}  \dV  &= \, \int_{\surf} f \, \extendDomain{\q} \dS  \formComma \quad \\
\lim\limits_{\varepsilon \to 0}\int_{\domain } \doubleWell(\phi) \partial_i \extendDomain{f} \, \extendDomain{\q}  \dV &= \,\int_{\surf}  \left( \Grad f(\Xb) \, \right)_i \extendDomain{\q}  \dS \formComma \quad \\
\lim\limits_{\varepsilon \to 0}\int_{\domain } \nabla \cdot \left( \doubleWell(\phi) \nabla \extendDomain{f}\, \right) \extendDomain{\q} \dV &=  \,\int_{\surf}  \Div \left( \Grad f \right) \extendDomain{\q}  \dS\formComma
\end{align*}
for $\extendDomain{\q} \in \Hi{\domain}$, see \cite{Raetz2007}.
A regularization is added to the function $\doubleWell(\phi)$ in some of the terms, to allow for a more stable solution of the linear system:
$\doubleWellRegularized(\phi) := \max(\doubleWell(\phi), \regularization)$ with $\regularization\ll 1$. This regularization is justified in \cite{Raetz2006, Li09}.

Applying a standard finite element method with globally continuous, piecewise linear elements $\mathbb{V}_h(\domain_h) = \{ v_h \in C^0(\domain_h) \; : \; v_h|_T \in \mathbb{P}^1, \, \forall \, T \in \domain_h \}$ on a triangulation $\Omega_h$ of $\Omega$, the time discretization as above and inserting the regularized delta function approximation $\doubleWellRegularized$, results in a sequence of diffuse interface problems: For $k=0,1,\ldots$, find $\extendDomain{\p}_{i}^{k+1}\in\mathbb{V}_h(\domain)$ s.t.
\begin{multline}
 \frac{1}{\tau_k} \int_{\domain_h} \doubleWellRegularized(\phi) \extendDomain{\p}_{i}^{k+1}\extendDomain{\q}\dV  \\
 + \K \int_{\domain_h} \doubleWellRegularized(\phi) \nabla \cdot \extendDomain{\pb}^{k+1} \partial_i \extendDomain{\q} + \doubleWell(\phi) \left[\nabla \cdot \left(\extendDomain{\pb}^{k+1} \times \nExt \right)\nabla \cdot \left( \vect{e}_i \extendDomain{\q} \times \nExt \right)  + \left(\extendDomain{\shapeOperator}^2 \cdot \extendDomain{\pb}\right)_i \extendDomain{\q}\right] \dV \\
 + \Kt\int_{\domain_h} \doubleWell(\phi) \nExt\cdot\extendDomain{\pb}^{k+1} \extendDomain{\surfNormalI}_{i} \extendDomain{\q} \dV + \Kn \int_{\domain_h} \doubleWell(\phi) \left[ \left( \|\extendDomain{\pb}^{k}\|^2 - 1 \right) \extendDomain{\p}_{i}^{k+1} + 2 \extendDomain{\p}_{i}^{k} \extendDomain{\pb}^{k}\cdot \extendDomain{\pb}^{k+1} \right] \extendDomain{\q} \dV \\
 = \frac{1}{\tau_k} \int_{\domain_h} \doubleWellRegularized(\phi) \extendDomain{\p}_{i}^{k}\extendDomain{\q}\dV + 2\Kn \int_{\domain_h} \doubleWell(\phi) \|\extendDomain{\pb}^{k}\|^2 \extendDomain{\p}_{i}^{k} \extendDomain{\q} \dV,  \quad \forall \; \extendDomain{\q} \in \mathbb{V}_h(\domain)\formComma
\end{multline}
for $i = 1, \ldots, 3$, with ${\extendDomain{\pb}}^0$ a smooth extension\footnote{A smooth extension to the domain $\Omega$ is implemented by successively extending fields to its surroundings, utilizing \eqref{eq:smooth_extension}, until the whole domain is covered, see also \cite{Stoecker2008}.} of $\pb^0$ to the domain $\Omega$. To assemble and solve the resulting system we use the FEM-toolbox AMDiS \cite{AMDiS2007,AMDiS2015} with domain decomposition on 64 processors. As linear solver we have used a restarted GMRES method with a restart cycle of 30, modified Gram-Schmidt orthogonalization, and a block Jacobi preconditioner with ILU($0$) local solver on each partition, as above for the \SFEM\ method.

\subsection{\label{sec:surface-approximation}Surface approximation and grids}
Surfaces similar to a sphere $\Sp$ can be triangulated by projecting a triangulation of the sphere to $\surf$, utilizing the coordinate projection $\pi$. For the \DEC\ method this triangulation must be well-centered, in other words, the circumcenter of each surface triangle must be located within the triangle. This property can be realized by triangles with internal angle less than $90^\circ$. An iterative procedure is applied to the projected sphere triangulation to fulfill this requirement, by shifting points tangentially to the surface so that all triangles have nearly equal internal angles and edge lengths. The algorithm is described in \cite{Nitschke2014}.

Other surfaces may be triangulated by cutting tetrahedra at the zero-level set of an implicit surface description. This triangulation must be optimized by retriangulation, \eg, by using \cite{ACVD, Valette2008}, and utilizing additionally the iterative procedure to get a well-centered complex, as above. Recently, an algorithm for mesh optimization, based on an edge collapsing strategy, was implemented in \cite{meshconv}. Even if \SFEM\ would need less requirements on the surface mesh we use the same meshes as for \DEC.
We have chosen a grid width $h$, \ie, the maximal edge length radius of all triangles, to be approximately $1/6$ of the defect core radius that is estimated experimentally.

For \DI\ we use a 3D conformal tetrahedral mesh adaptively refined near the surface. Therefore, the interfacial region, \ie, $\{\Xb\in\Omega\,:\,\phase(\xb)\in[0.1, 0.9]\}$, contains approximately  $7$ grid points in normal direction to the surface. This refinement guarantees good agreement with the sharp surface limit, see, \eg, \cite{Aland2010, Article2013} for a justification and quantitative study. The signed-distance function, the phase-field is based on, is calculated from the triangulated surface by an algorithm utilizing a ray tracing principle. For every grid point in the 3D mesh the distance to the surface is calculated and afterwards the correct sign is assigned. This algorithm is explained and implemented in \cite{meshconv} and has an asymptotic complexity of $\sim\mathcal{O}(|\Omega_h|\cdot\log\numFs)$.

\section{Computational results}\label{sec:ComputationalResults}
We validate the proposed approaches on the unit sphere. Due to lack of analytical description of minimizers $\pb \in \Hdr{\surf}{\Tangent \surf}$, we compare the numerical results with each other. The \DEC\ approach thereby serves as reference. We also explore the stability of minimal energy defect configurations on more complicated surfaces with non-constant curvature and demonstrate the tight interplay of defect localization and geometric properties. Within these studies we show the possibility of equilibrium states other than the trivial realization of the Poincar\'e-Hopf theorem and thus the possibility to reduce the weak surface Frank-Oseen energy by incorporating additional defects. To validate these results we again compare the numerical results with each other. The penalty parameter $\Kn$ is chosen such that the defect core radius is resolved, see \autoref{tbl:parameters}. The section is concluded by providing information on the numerical effort for each method.

\begin{table}
\begin{center}
\begin{tabular}{|r|l|l|l|}
\hline
& & sphere & nonic surface \\\hline
time & $t_\text{end}$ & $5$   & variable ${(\ast)}$ \\
     & $\tau_k\equiv\tau$ & $10^{-3}$ & $5\cdot 10^{-4}$  \\ %
model & $\K$  & $1$           & $1$ \\
      & $\Kn$ & $10^3$        & $200$ \\
      & $\Kt$ & $[10^3-10^5]$ & $10^5$ \\\hline
\SPH  & $N$   & $190$         & --- \\
      & $N_{\lu}$ & $250$     & --- \\
      & $N_{\lv}$ & $400$     & --- \\
      & $\tau$ & $2\cdot 10^{-4}$     & --- \\ %
\SFEM & $h$       & $0.013$ & $0.035$ \\
\DI   & $\varepsilon$ & $0.15$ & $0.2$ \\
      & $\regularization$ & $10^{-6}$  & $10^{-6}$ \\
      & $h$       & $0.023$ & $0.078$ \\
      & $\Omega$  & $[-1.5, 1.5]^3$ & $[-2,3]\times[-2,2]^2$ \\\hline
\end{tabular}
\end{center}
\caption{\label{tbl:parameters}Simulation parameters for the two setups: relaxation on the sphere and nonic surface. $(\ast)$ the end-time of the simulation is chosen so that the system is close to equilibrium, \ie, if the criterion $|\FF(t_{k+1}) - \FF(t_k)| < 10^{-14}\cdot|\FF(t_{k+1})|$ is fulfilled.}
\end{table}

\subsection{Method comparison on sphere}\label{sec:ModelComparison}
We consider an initial condition $\pb^0$ with two sinks (\includegraphics[height=8pt]{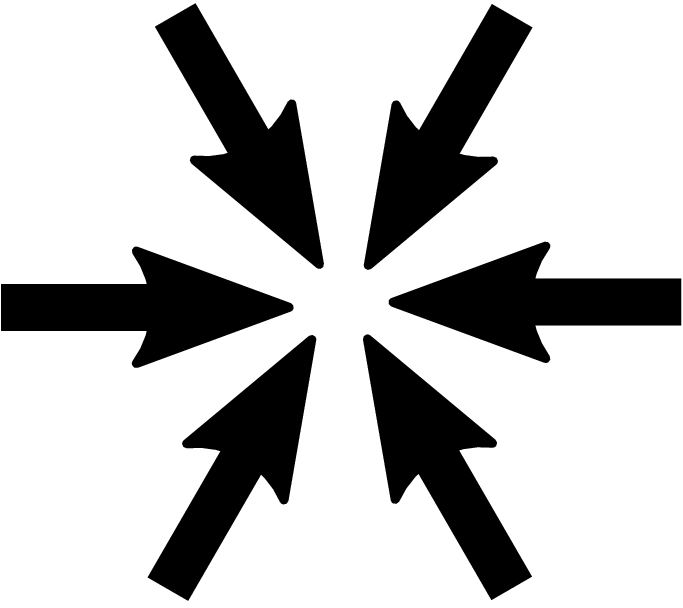} $+1$), a source (\includegraphics[height=8pt]{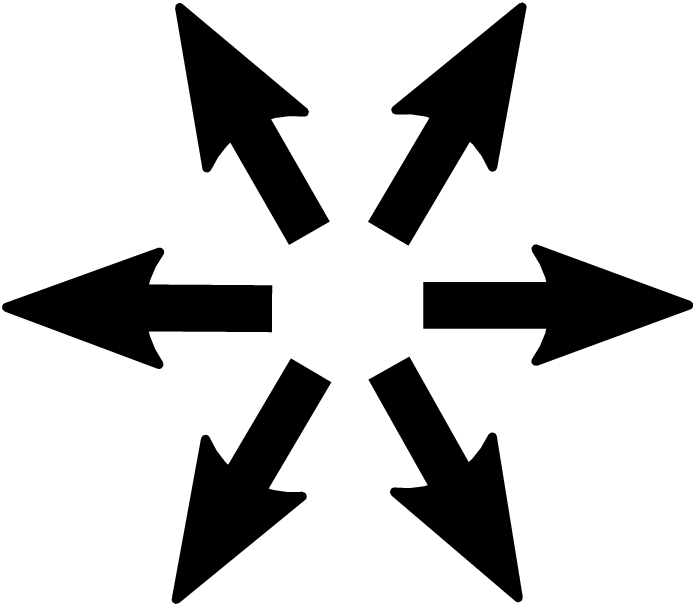} $+1$) and a saddle point (\includegraphics[height=8pt]{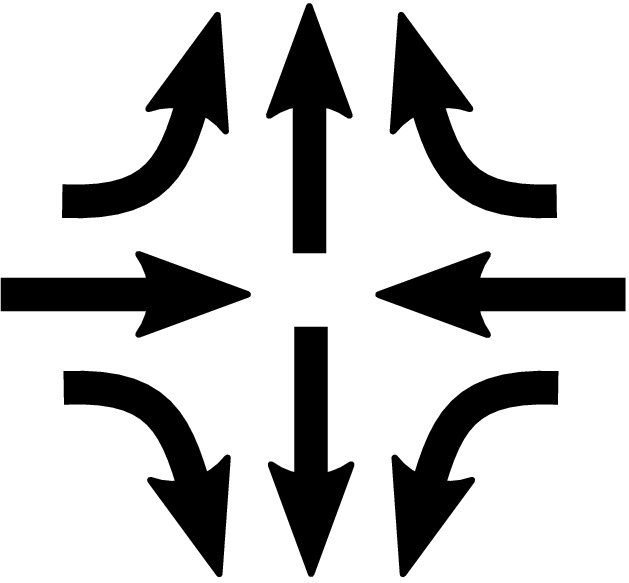} $-1$) on the unit sphere $\surf=\Sp$. The numbers are the topological charges or the winding numbers $\text{ind}_V(\mathbf{d}_i)$ of the defects $\mathbf{d}_i$. They are defined as the
algebraic sum of the number of revolution of $\pb$ along a small counterclockwise oriented curve around the defect.
The Poincar\'e-Hopf theorem requires
\begin{align}
\sum_i \text{ind}_V(\mathbf{d}_i) &= \chi(\surf),
\end{align}
which in the present case is satisfied as $1 + 1 + 1 - 1 = 2$. The four defects are positioned equidistant on the
$x$-$y$-equatorial plane. To avoid metastable configurations we shift one sink defect
slightly closer to the saddle point defect.
\begin{align}\label{eq:InitCondition}
\begin{aligned}
  \pb^0 &= \frac{\ProjectSurf \hat{\pb}^0}{\|\ProjectSurf \hat{\pb}^0\|}\formComma \quad \mbox{ where }\\
\hat{\pb}^0 &=
\begin{cases}
  \left[-x , 0, -z\right]^T  & |y| \geq \cos\frac{\pi}{4}  \\
  \left[ 0 , y , z\right]^T  & x \geq \cos\frac{\pi}{4}\\
 \left[0 , \sin\left(\pi \left( y - \lambda\right)\right), -\sin\left( \pi z \right) \right]^T  & x \leq -\cos\frac{\pi}{4}\\
  \left[ \left|\frac{y}{\cos\frac{\pi}{4}}\right| - 1, \frac{y}{\cos\frac{\pi}{4}}, 0\right]^T & \mbox{otherwise}
\end{cases}
\end{aligned}
\end{align}
with $\lambda = 0.01$ used in our simulations.

\begin{figure}[ht]
 \begin{center}
 \includegraphics{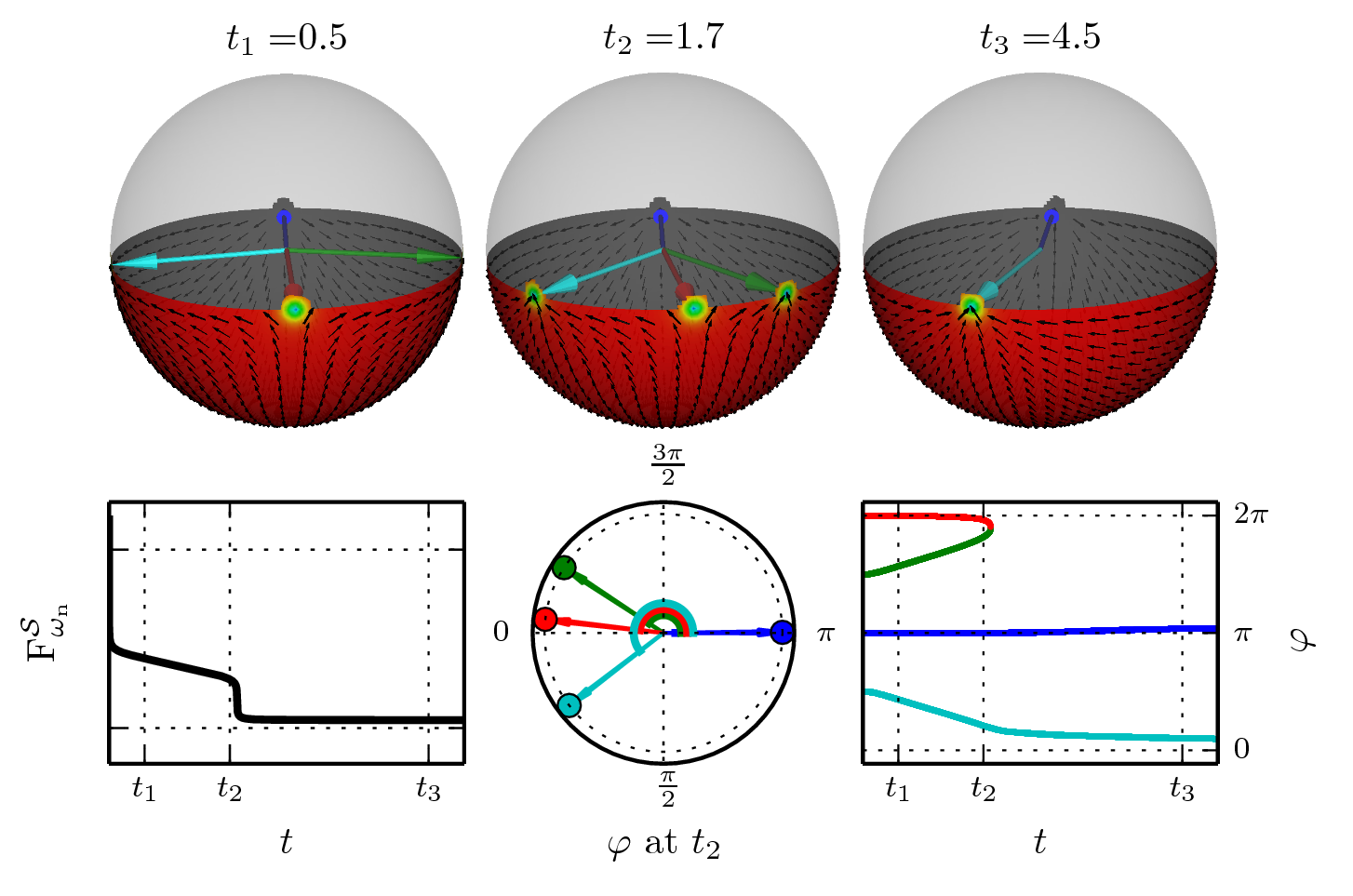}
  \caption{\label{fig:FourDefectsTrajectories}(Colors online) Top: Sequence of director field configurations (glyphs) and defect positions (color gradient on surface and large arrows) in
  the evolution of the four-defect test case, at time $t_0$, the four defect configuration,
  time $t_1$, the defect annihilation, and time $t_2$, the two-defect configuration. Bottom:
  Energy evolution (left), defect positions in $x$-$y$-equatorial plane at $t_2$ (middle) and defect trajectories of the four-defect test case (right).
  The angle $\lv \in [0,2\pi)$ describes the defect positions in the $x$-$y$-equatorial
  plane. Colors of the defects: source (dark blue), sinks (cyan and green), saddle point (red).}
 \end{center}
\end{figure}

Since opposing topological charges attract each other we observe the motion of
the two sink defects to the saddle point defect and eventually the fusion of
the saddle point defect with the closer sink defect (see \autoref{fig:FourDefectsTrajectories}).
The time needed for the annihilation of the two defects is denoted by $t_f$ and called fusion time.
Finally, the remaining two defects relax to a position with maximal distance.
Due to the symmetry of the setup the defect positions will remain in the equatorial
plane.

These dynamics are consistently observed within all methods. To measure deviations in the proposed numerical methods we compare against the
\DEC\ solution. Therefore, we introduce as quantitative measure a density like mean energy error $\epsilon_\text{e}$ (normalized by the area $A$ of the surface, $A_{\Sp}=4\pi$) and
as qualitative measure the error in the defect fusion time $\epsilon_\text{f}$,
\begin{align}
 \label{eq:errMeassure}
 \epsilon_\text{e} &:= \frac{1}{A\, t_\text{end}}\int_{0}^{t_\text{end}} \left| \frac{\F{\Kn,(\Kt)}^{\surf (\Omega)}(\text{\shortcut{M}}) - \F{\Kn}^\surf(\text{\DEC})}{\F{\Kn}^\surf(\text{\DEC})}\right| \dif t \formComma \\
 \epsilon_\text{f} &:= \left|\frac{t_f(\text{\shortcut{M}}) - t_f(\text{\DEC})}{t_f(\text{\DEC})}\right| \formComma
\end{align}
for a numerical method \shortcut{M}. Within this framework we evaluate the proposed vector-valued methods \DEC\ and \SPH, and the component-wise methods \SFEM\ and \DI, with parameters from \autoref{tbl:parameters}.

\begin{figure} [ht]
 \begin{center}
 \includegraphics{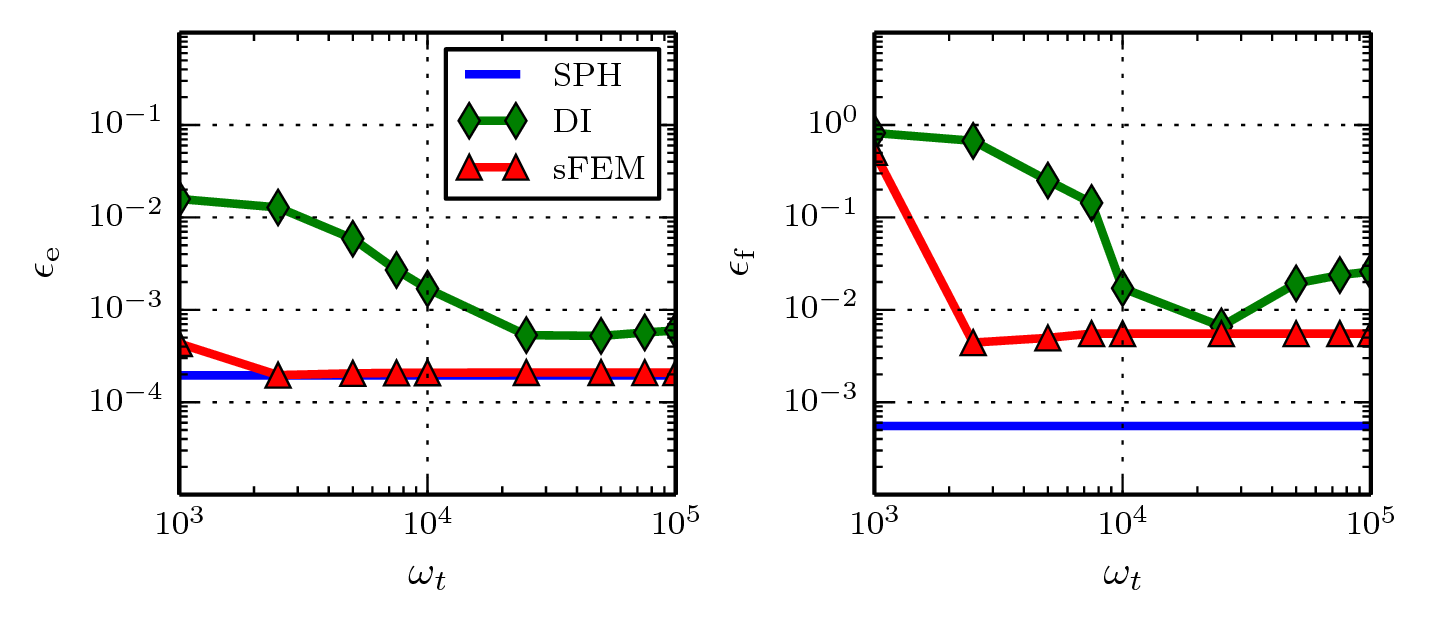}
  \caption{\label{fig:RelErrorQuantQualiSphere}
 (Colors online) The errors in mean energy $\epsilon_\text{e}$ (left) and fusion time of
  defects $\epsilon_\text{f}$ (right) for various tangentiality penalty parameters $\Kt$.
  Three different methods are compared to DEC: spherical harmonics (blue), parametric FEM (red),
  and diffuse interface (green).}
 \end{center}
\end{figure}

\autoref{fig:RelErrorQuantQualiSphere} shows the obtained computed errors. The methods essentially show matching solutions. The relative energy difference and difference in defect fusion time is reduced for increasing penalty factor $\Kt$, but is limited by the differences in the compared methods, \eg, difference in the location of DOFs and the discretization of the surface.
The \SPH\ method does not depend on a tangentiality penalization as the \DEC\ method. Thus, the error values result from a difference in the surface representation and the truncation in the spherical harmonics expansion.
Apart from this, two qualitatively different behaviors for \SFEM\ and \DI\ can be observed. Where the method \SFEM\ shows nearly constant errors (at least for $\Kt > 2500$), the method \DI\ shows a dependence on the penalty parameter.
This effect arises from the interaction of the penalty forcing and the geometric approximation of $\surf$ by a smeared-out delta-function, \ie, a non-constant penalty factor throughout the interface. Close to the surface the director field $\extendDomain{\pb}$ is not guaranteed to be tangential to $\surf$ for $\Kt$ too small. Increasing the penalty factor finally leads to tangential fields in the surrounding of the interface. This results in error values close to those of \SFEM. A difference in these two methods is expected, due to the additional approximation of the surface and the surface differential operators by the diffuse interface representation.

Within a reasonable tolerance all four methods show the same dynamic behavior along quantitative and qualitative computed errors and converge to the same stationary solution with two defects, a source (+1) and a sink (+1), which are at maximal distance from each other.
 
\subsection{Higher order surfaces}\label{sec:HigherOrderSurfaces}
\begin{figure} [ht]
 \begin{center}
  \includegraphics{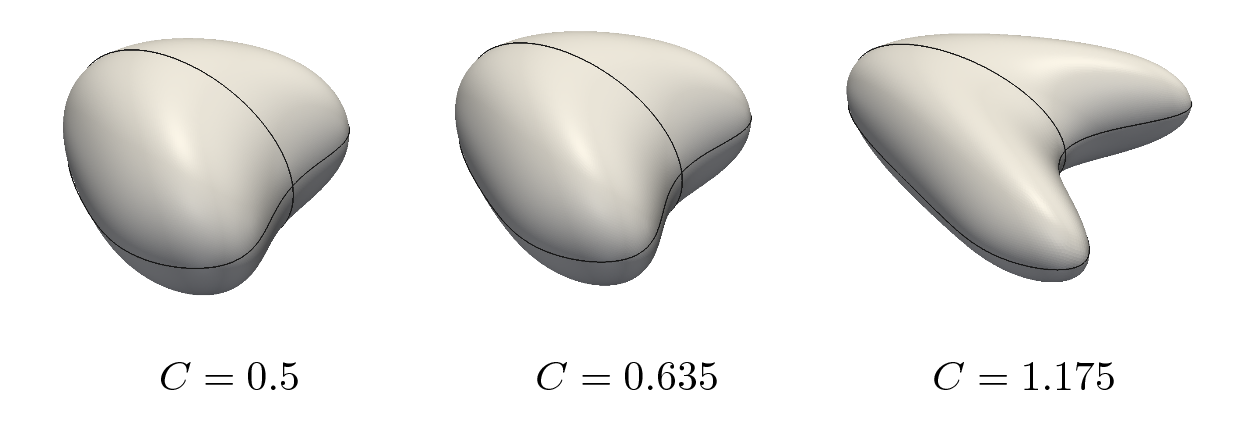}
  \caption{\label{fig:nonic_surfaces}Nonic surfaces corresponding to three different
  stretching parameters $C$. Left: surface with defect fusion-time $> 0$, Center:
  four-defect configuration gets stable, Right: four-defect configuration is
  energetically equivalent to two-defect configuration.}
 \end{center}
\end{figure}

To further validate the consistency of the methods \DEC, \SFEM\, and \DI, we extend the test-setup to a sequence of surfaces with non-constant
curvature, see \autoref{fig:nonic_surfaces} for examples. All surfaces have $\chi(\surf) = 2$, thus allowing defect configurations as in the previous example.

The construction of the surfaces is based on a deformation of the unit sphere, such
that regions with positive and negative Gaussian curvature emerge. Our goal is to study the influence of these regions on the defect location. Are defects attracted by
these regions? Is there a relation between the topological charge of the defect and the Gaussian curvature?

The postulated parametrization of the unit sphere \( \Sp \), $\Xb_{\Sp}(\lu,\lv)$, given in \autoref{sec:SpectralMethod}, is stretched in the
\( z \)-direction by the displacement function
$f_{\pstretch,\pprop}$ with factors \( \pprop\in (0,1) \) and \( \pstretch > 0 \),
\begin{align*}
  f_{\pstretch,\pprop}(z) &:= \frac{1}{4} \pstretch z^2 \left[(z+1)^2 (4-3 z)+\pprop (z-1)^2 (4+3 z)\right]
\end{align*}
and compressed along the \( y \)-direction by a factor \( \ppress\in [0,1) \). This leads to the parametrization
\begin{align}
  \Xb(\lu,\lv) &:= \Xb_{\Sp}(\lu,\lv) + f_{\pstretch,\pprop}(\cos\lu)\EuBase{x} - B\sin\lu\sin\lv\EuBase{y}
  \formPeriod
\end{align}
The surface can also be expressed implicitly by the zero-level set of the function
\begin{align}\label{eq:nonicImplicite}
  \varrho(x,y,z) &:= \left( x - f_{\pstretch,\pprop}(z) \right)^{2} + \frac{1}{(1-\ppress)^{2}}y^{2} + z^{2} - 1\formPeriod
\end{align}
This gives a polynomial $\varrho$ of degree 10, which motivates the name \emph{nonic surfaces}. The asymmetry of the surfaces \wrt\ the \( x \)-\( z \)-plane prevents metastable defects configurations. The necessary surface quantities can be derived directly from the level set formulation by $\extendDomain{\surfNormal} = \nabla \varrho / \|\nabla\varrho\|$ and $\shapeOperator_{ij} = -\left[\ProjectSurf \nabla \extendDomain{\surfNormal}_j\right]_i $.

To investigate the energy value $\E{\pb^\ast}$ of a stationary solution $\pb^\ast$
and the stability of defect configurations we analyze the evolution of two different
initial solutions $\pb^{0}_{(4)}$ and $\pb^{0}_{(2)}$. The first one, $\pb^{0}_{(4)}$, has four separated defects, while the second one, $\pb^{0}_{(2)}$, has two.

\begin{figure} [ht]
 \begin{center}
  \includegraphics{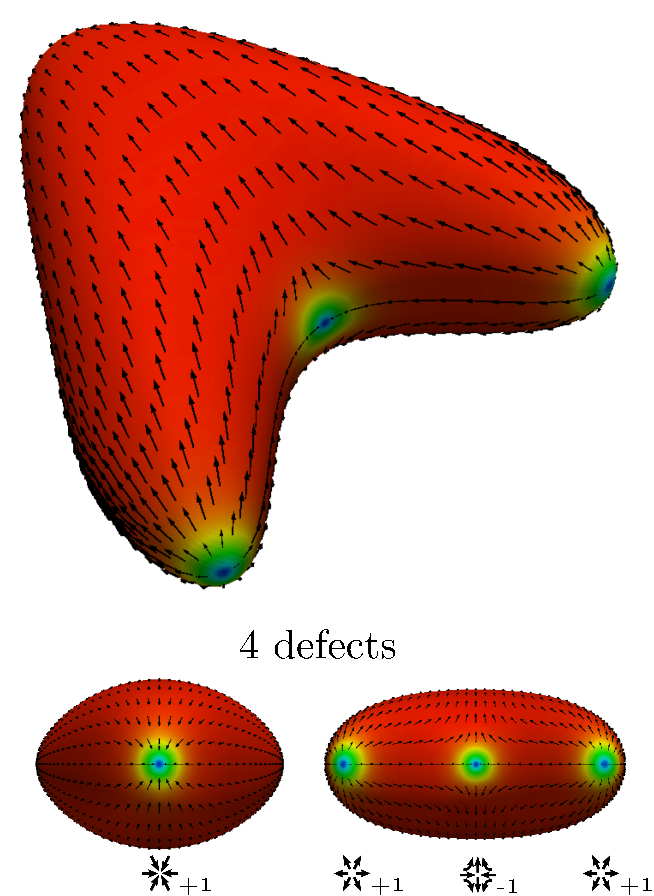}
  \includegraphics{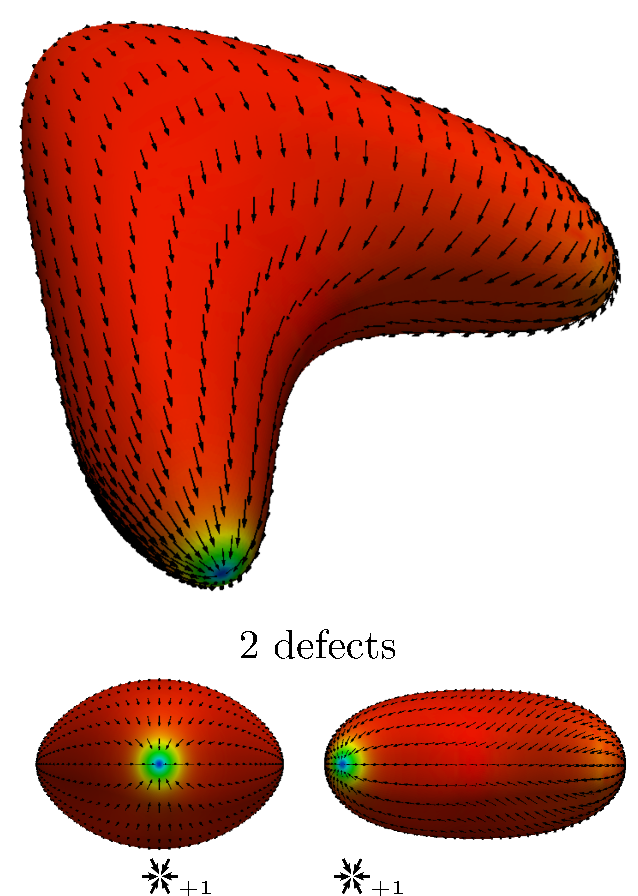}
  \caption{\label{fig:nonicDefects}Equilibrium states for surface with $\pstretch=1.175$, norm defects (color gradient) and director (glyphs). Second row: back and front detail of configuration.}
 \end{center}
\end{figure}

At first, we consider the projected unit vector $\EuBase{x}$, which can be represented by the
surface gradient of the \( x \)-coordinate, \ie,
\begin{align}
  \pb^{0}_{(4)} := \ProjectSurf\EuBase{x} = \Grad x = (\exd x)^{\sharp} \formPeriod
\end{align}
On an edge \( e=\left[ v_{1} , v_{2} \right]\in\Es \),
where the face \( \face_{1}\succ e \) is right of \( e \) and \( \face_{2}\succ e \) is left of \( e \),
so that \( \star e = [c(\face_{1}),c(e)] + [c(e),c(\face_{2})] \) is the dual edge,
we can approximate the 1-form \( \exd x \), utilizing integration by parts on  \( e \), by
\begin{align}
  \alphav^{0}_{(4)}(e) &=
      \left( v_{2}^{x} - v_{1}^{x},
	  -\frac{\left| e \right|}{\left| \star e \right|} \left( \left[c(\face_{2})  \right]^{x} - \left[c(\face_{1})  \right]^{x} \right)  \right)
  \formPeriod
\end{align}

To enforce a two-defect solution in equilibrium for the second case, we project a
slightly rotated unit vector \( \EuBase{y} \) to the surface. The rotation by an
angle \( \gamma \) in the normal plane of the \( \R^{3} \)-vector \( [-1,0,1]^{T} \)
is thereby represented by the rotation matrix \( R_{\gamma} \). This defines
\begin{align}
  \pb^{0}_{(2)} &:= \ProjectSurf R_{\gamma} \EuBase{y} \formPeriod
\end{align}
Our choice of \( \gamma \) is \( 0.05 \). In the context of \DEC, the evaluation
of a vector field \( \qb\in\Tangent\surf \) with the dual edge vector
\( \eb_{\star} \) on edge \( e \) at the intersection \( e\cap\star e = c(e) \)
is ambiguous. To overcome this, we define in a canonical way a dual 1-chain, utilizing
the definition of a dual edge \( \star e = \star e|_{\face_{1}} + \star e|_{\face_{2}}\).
This leads to
\begin{align*}
  \qb(c(e))\cdot\eb_{\star} :=  \qb(c(e)) \cdot \left( \eb_{\star}|_{\face_{1}} +  \eb_{\star}|_{\face_{2}}\right)
            = \qb(c(e)) \cdot \left( c(\face_{2}) - c(\face_{1}) \right) \formComma
\end{align*}
where the face \( \face_{1}\succ e \) is right of the edge \( e \) and \( \face_{2}\succ e \) is located left.
Thus we get the initial discrete PD-1-form
\begin{align}
  \alphav^{0}_{(2)}(e) = \left( \pb^{0}_{(2)}(c(e)) \cdot \eb ,
	    -\frac{\left| e \right|}{\left| \star e \right|}\pb^{0}_{(2)}(c(e)) \cdot \eb_{\star} \right)
      \formPeriod
\end{align}
The normalized versions of $\pb^{0}_{(i)}$ and $\alphav^{0}_{(i)}$ can easily be
constructed by point-wise or edge-wise normalization, respectively, using the
definition of the norm in \eqref{eq:PDNorm} for the discrete PD-1-forms.

Within this setup we evaluate the energy for stationary solutions  $\pb^\ast$ and the
number of defects for both initial solutions $\pb^{0}_{(4/2)}$ for a sequence of values
$\pstretch \in [0,1.5]$. The parameter $\pprop=0.95$ remains fixed while $\ppress$ is
related to $C$ by $\ppress = 7/20\pstretch$.

An example of the two different initial fields relaxed to equilibrium is shown in \autoref{fig:nonicDefects} for a specific nonic surface.
We find $+1$ defects at extrema of the Gaussian curvature, while a $-1$ defect may appear at the saddle point.
This dependency is in agreement with results for the similar problem of flow on curved surfaces \cite{Reuther2015,Nitschke2016}.

\begin{figure} [ht]
 \begin{center}
  \includegraphics{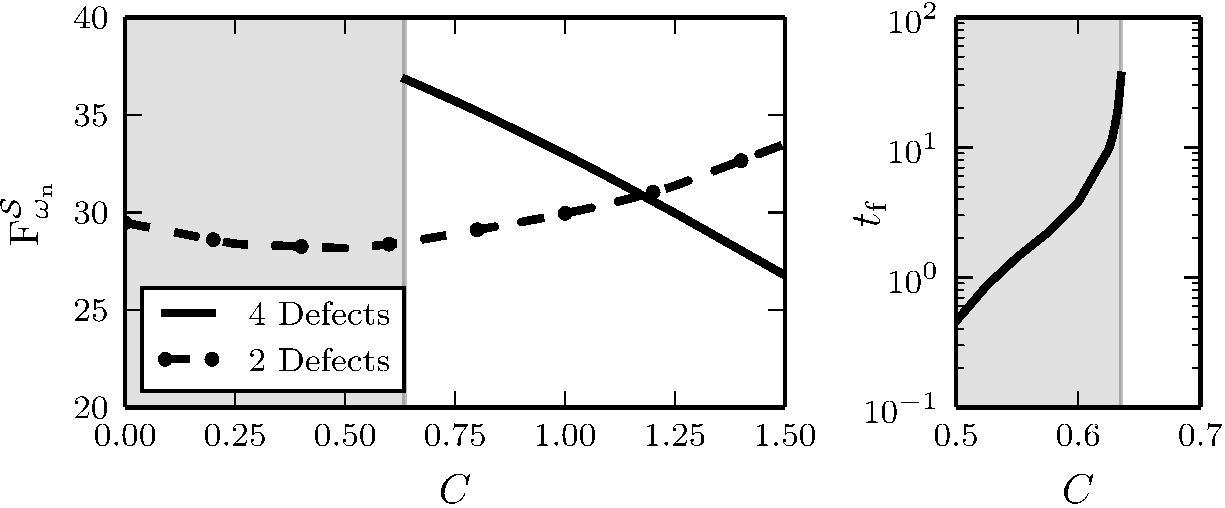}
  \caption{\label{fig:nonicEnergies}Energy $\EE$ for stationary solutions with four and
  two defects for nonic shapes with $\pstretch \in [0,1.5]$ (left) and defect fusion time
  for the four-defect initial solution (right).}
 \end{center}
\end{figure}

For shapes with $\pstretch \in [0.5, 0.635]$ we observe that both initial solutions
converge to a two-defect configuration. In \autoref{fig:nonicEnergies} (right)
we plot the fusion time for defect annihilation for initial condition $\pb^{0}_{(4)}$.
Notice the steep increase in this time for $\pstretch \nearrow 0.635$.
For $\pstretch \gtrsim 0.635$ a four-defect configuration becomes stable. It poses a local energetic minimum. Further increasing the parameter $\pstretch$, continuously amplifies
the Gaussian curvature on the bulges and saddle. As shown in \autoref{fig:nonicEnergies} (left), this leads to a decreasing energy cost for the four-defect stationary solution, while costs for the two-defect solution increase monotonically until the energies are equal at $\pstretch \approx 1.175$. For $\pstretch \gtrsim 1.175$ the four-defect solution
becomes energetically favorable. This behavior is stable against variations in the penalty parameter $\Kn$, which is chosen, such that the defect core radius is resolved, see \autoref{tbl:parameters}.

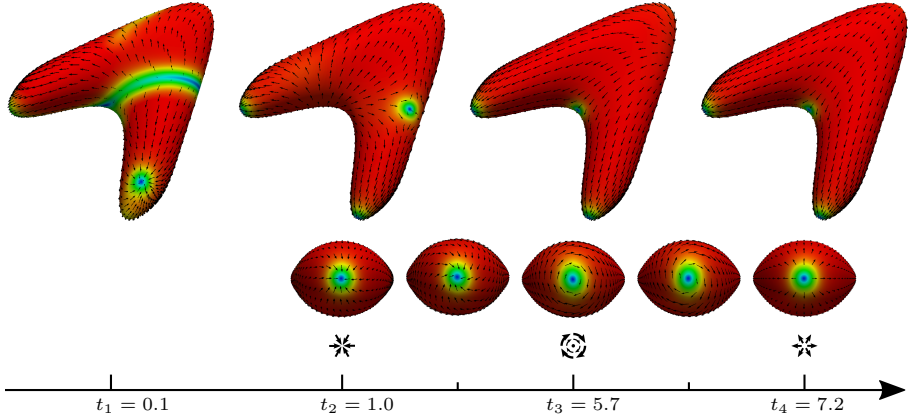
\begin{figure} [ht]
 \begin{center}
\begingroup%
  \makeatletter%
  \providecommand\color[2][]{%
    \errmessage{(Inkscape) Color is used for the text in Inkscape, but the package 'color.sty' is not loaded}%
    \renewcommand\color[2][]{}%
  }%
  \providecommand\transparent[1]{%
    \errmessage{(Inkscape) Transparency is used (non-zero) for the text in Inkscape, but the package 'transparent.sty' is not loaded}%
    \renewcommand\transparent[1]{}%
  }%
  \providecommand\rotatebox[2]{#2}%
  \ifx\svgwidth\undefined%
    \setlength{\unitlength}{12cm}%
    \ifx\svgscale\undefined%
      \relax%
    \else%
      \setlength{\unitlength}{\unitlength * \real{\svgscale}}%
    \fi%
  \else%
    \setlength{\unitlength}{\svgwidth}%
  \fi%
  \global\let\svgwidth\undefined%
  \global\let\svgscale\undefined%
  \makeatother%
  \begin{picture}(1,0.45)%
    \put(0,0){%
\includegraphics[width=\unitlength]{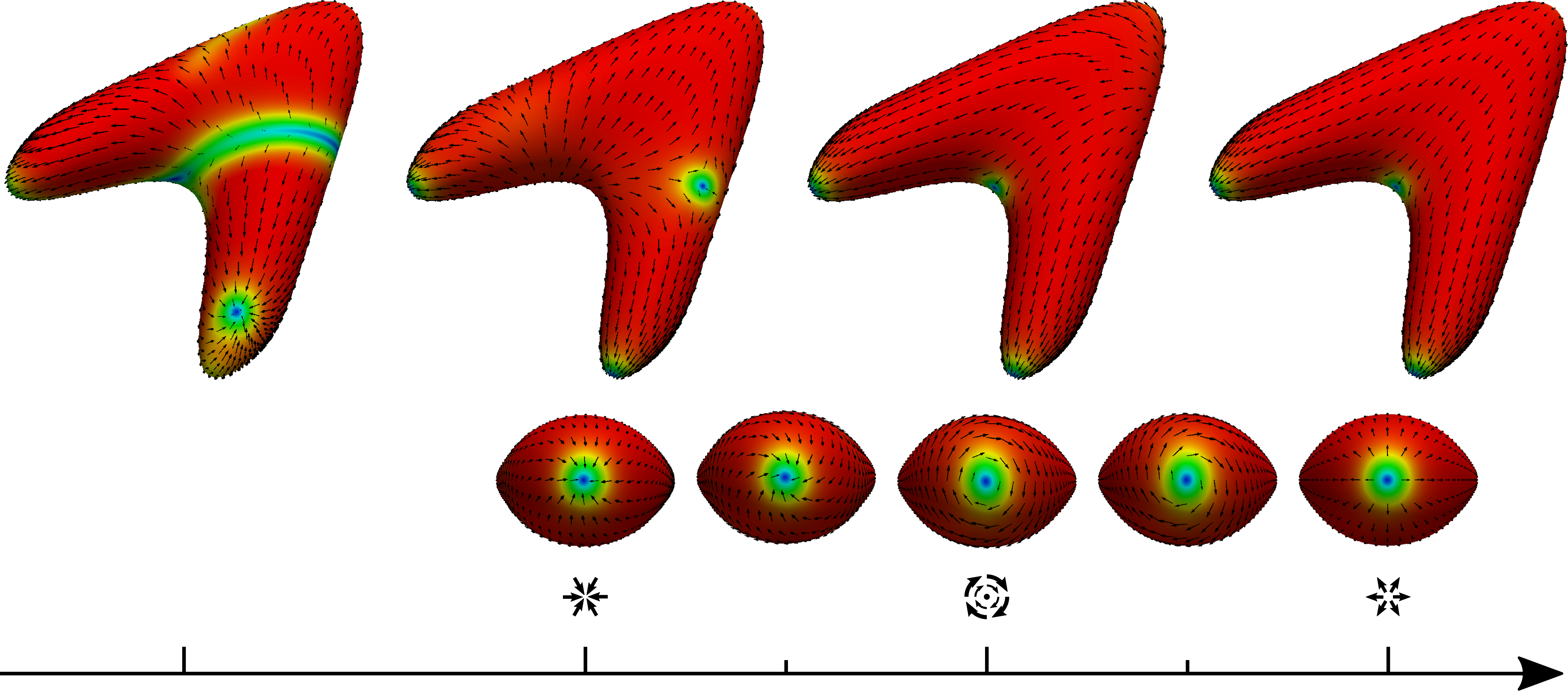}%
    }%
    \put(0.10029845,-0.015){\color[rgb]{0,0,0}\makebox(0,0)[lb]{\scriptsize{$t_1 = 0.1$}}}%
    \put(0.34889693,-0.015){\color[rgb]{0,0,0}\makebox(0,0)[lb]{\scriptsize{$t_2 = 1.0$}}}%
    \put(0.59914185,-0.015){\color[rgb]{0,0,0}\makebox(0,0)[lb]{\scriptsize{$t_3 = 5.7$}}}%
    \put(0.84760166,-0.015){\color[rgb]{0,0,0}\makebox(0,0)[lb]{\scriptsize{$t_4 = 7.2$}}}%
  \end{picture}%
\endgroup%
   \caption{\label{fig:nonicDynamics}Snapshots of the time evolution and the final stationary solutions with four defects on the nonic shapes with $\pstretch = 1.5$. First row: top view. Second row: back view with a single defect, evolution from a sink shape (at $t=1$) over vortex shape (at $t=5.7$) to the final source shape (at $t=7.2$).}
 \end{center}
\end{figure}

These experiments emphasize the impact of curvature on the energetic cost of a defect configuration and prove the key role of domain geometry in enabling non-trivial realizations of the
Poincar\'e-Hopf theorem. \autoref{fig:nonicDynamics} shows snapshots of the evolution on the most deformed surface with $C = 1.5$ and noise used as initial condition. Which stationary shape is selected strongly depends on the initial condition. We here only show the one converging to the four-defect configuration.

\begin{figure} [ht]
 \begin{center}
  \includegraphics[scale=0.8]{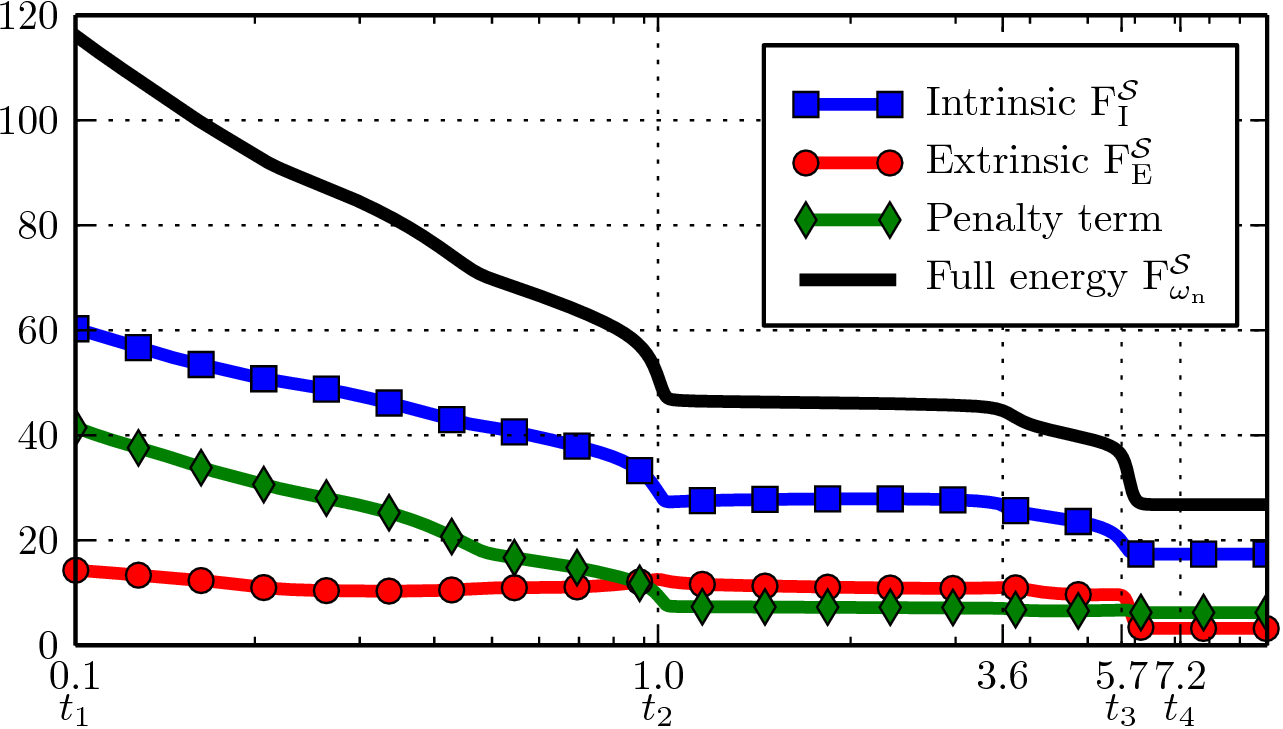}
  \caption{\label{fig:nonicDynamicsEnergy}(Colors online) Development of the energy parts in the relaxation starting from random initial state. The four parts plotted with lines and symbols sum up to the full energy $\EE$. Highlighted are five time steps that mark changing events. From $t=0.1$ to $t=3.6$ the defects move to their final position. In the time period $t=1.0$ to $t=5.7$ the back defect rotates by 90 degree, from a sink to a vortex defect. From time $t=5.7$ to $t=7.2$ the back defect rotates further by 90 degree, from a vortex defect to a source defect.}
 \end{center}
\end{figure}

The relaxation shows four periods with distinct behavior. Starting from a random initial configuration the noise smoothes out to a state with emerging localized defects at time around $t=0.1$. Until time $t=1$ these defects reach their final normalization shape, in other words, the penalization term in the energy reduces up to this time and stays constant from this time on, as can be seen in \autoref{fig:nonicDynamicsEnergy}. The defects move at first slowly and then very fast to their final position around the high curvature areas and the saddle point. This happens until time $t=3.6$. When the back defect reaches its final position it starts to rotate the vector field up to 90 degrees. Thus, a sink defect evolves to a vortex defect at around time $t=5.7$. This process continues and rotates the vector field around this back defect further by 90 degrees until a source defect shape is reached at around time $t=7.2$.

\begin{figure} [ht]
 \begin{center}
  \includegraphics{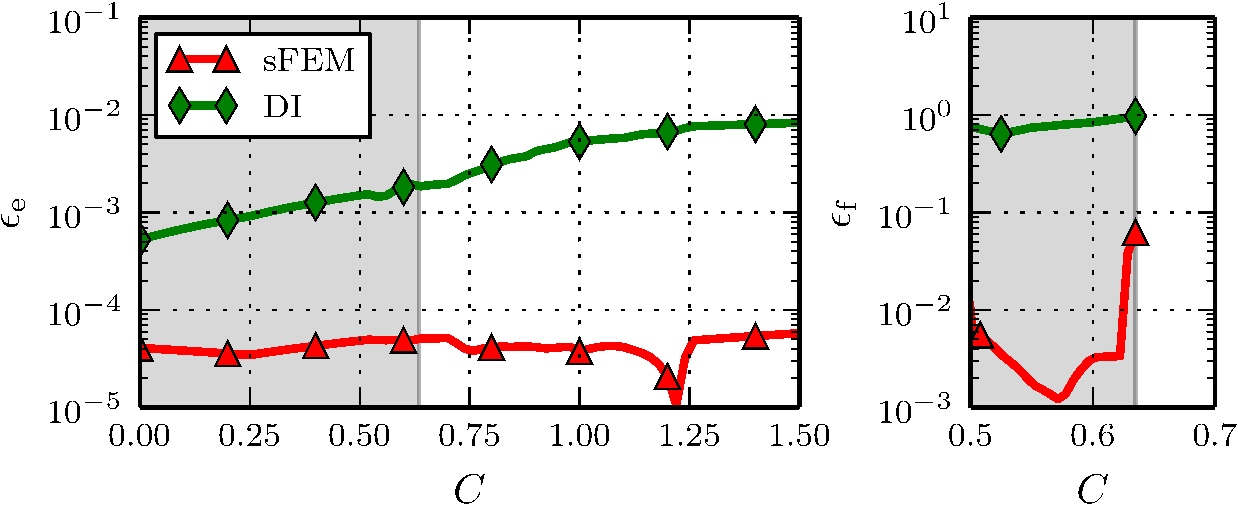}
  \caption{\label{fig:nonicErrors}(Colors online) Relative errors \wrt\ DEC solution
  of mean energy (left) and fusion time (right) for nonic shapes $\pstretch \in [0,2]$
  and the numerically methods parametric FEM (red dashed) and diffuse interface (green solid).}
 \end{center}
\end{figure}

Beside these exploratory results, shown in \autoref{fig:nonicEnergies} we also use this parameter study to verify the quality of the numerical methods \SFEM\ and \DI. In \autoref{fig:nonicErrors} we plot the relative errors introduced in \eqref{eq:errMeassure} for the mean energy and fusion time. As numerical parameters we
have chosen values listed in \autoref{tbl:parameters} in the column \textit{nonic surface}.

As shown in \autoref{fig:nonicErrors} (left) we observe the same behavior with both methods, across the full range of shapes $\pstretch \in [0, 1.5]$, within reasonable error bounds.
The more approximative \DI\ yields significant stronger deviations from the
\DEC\ results, up to two orders of magnitude in the mean energy error. Furthermore, we notice
increasing errors with amplified curvature. The critical point
$\pstretch_{\textup{crit}} = 0.635$ of emergence of a new stable defect configuration is
qualitatively reproduced by both methods. \DEC\ and \SFEM\ yield identical results for $\pstretch_{\textup{crit}}$, up to the
probing grid spacing of $\delta \pstretch = 2.5\cdot 10^{-3}$. \DI\ produces a
critical value of $0.7125$, which corresponds to a relative error of $0.122$ \wrt\
the \DEC\ result. As a result the dynamics evaluated by \DI\ close to this critical event
exhibit distinct deviations leading to substantial relative errors for the fusion
time as shown in \autoref{fig:nonicErrors} (right).

We do not compare the dynamic evolution if started from noise, as identical initial conditions cannot be specified. However,
also \SFEM\ and \DI\ produce evolutions which are qualitatively the same as in \autoref{fig:nonicDynamics}. Again, whether a two-defect or four-defect configuration is reached strongly depends on the initial condition.
 
\subsection{Performance comparison}\label{sec:PerformanceComparison}
We summarize pros and cons of the considered numerical methods, with respect to complexity, accuracy, generality and numerical performance.

As a first quality measure, we consider the applicability of the methods to various geometric surfaces. Here \SPH\ is the most restrictive as it can only be applied to spherical surfaces, since eigenfunctions and eigenvalues of the \LaplaceDeRham operator are utilized. \DEC\ and \SFEM\ can be applied to all surfaces, where a suitable surface mesh is available. \DEC\ requires well-centered simplicial surface elements, whereas the requirements for \SFEM\ are less restrictive. However, a non-regular shape of the triangles may increase the condition number of the resulting linear system \cite{Olshanskii2013,Dziuk2013}. Thus, the quality of the surface triangulation matters for both approaches. \DI\ uses an implicit description of the surface and thus does not rely on an approximate surface mesh. The 3D domain $\Omega$ can be adaptively triangulated using regular shaped tetrahedra and thus allows to conserve good mesh quality easily. Efficient methods to calculate a signed-distance function $d_\surf$ from an implicit description of $\surf$ or from a triangulated surface are necessary and available for tetrahedral meshes, see \cite{Bornemann2006,Stoecker2008}.

The computational costs for all the methods vary a lot. Denoting by $\numVs$ the number of vertices of a surface triangulation and by $\numEs$ the number of edges. For \SPH\ the main computational expenses are related to the forward and backward transform, which can be classified as $\mathcal{O}(N^2\log N + \numVs)$ with band-width $N$, typically $N\sim\sqrt{\numVs}$. The other methods have to assemble and invert a linear system in each time step iteration. The number of degrees of freedom (DOFs) and the corresponding average number of non-zero entries (NNZ) per row in the linear system are summarized in \autoref{tbl:complexity}. The total number of non-zeros in the system is approximately the same for \DEC\ and \SFEM, whereas \DI\ produces a much larger and denser system.

\begin{table}
\begin{center}
 \begin{tabular}{|c|c|c|c|} \hline
   & \DEC & \SFEM & \DI \\\hline
     $\sharp$DOFs & $2\cdot\numEs \approx 6\cdot\numVs$ & $3\cdot\numVs$ & $\gg 3\cdot\numVs$ \\
     NNZ/row & $12$ & $20$ & $37$ \\\hline
 \end{tabular}
 \caption{\label{tbl:complexity}Number of degrees of freedom $\sharp$DOFs and number of non-zeros per row of the matrix for the three methods that assemble a linear system.}%
\end{center}
\end{table}

The structure of the linear systems is also different. Where the \SFEM\ and \DI\ method produce symmetric matrices for symmetric differential operators, the \DEC\ approach results in a non-symmetric matrix, since not all triangles in the discretization are equilateral. This restricts the choice of linear solvers and often results in an additional performance overhead.

\DI\ allows to use classical finite element software. The additional cost, resulting from the treatment in 3D can be reduced by adaptive refinement in a narrow band around the surface. This establishes this approach as an easy to use tool also in the context of surface vector field calculations. A further extension of the analyzed models toward evolving surfaces can also most easily be adopted to \DI\ methods by evolving the implicit function or the phase-field variable.

\section{Conclusion and Outlook}\label{sec:ConclusionOutlook}
We presented a brief derivation of the weak surface Frank-Oseen energy as a thin film limit of the well known 3D Frank-Oseen distortion energy. By penalizing the unity of the vector field the limit can be established for surfaces with $\chi(\surf) \neq 0$. We highlight the importance of intrinsic and extrinsic energy contributions. Dynamic equations for surface bound polar order are obtained by an $L^2$-gradient flow approach, leading to a vector-valued surface PDE.

The energy and the dynamic equations have been adapted to suit several numerical methods. The least approximating methods base on a direct discretization of the vector-valued state space of the energy functional. For spherical surfaces this is \SPH\ and for arbitrary surfaces \DEC. Extending the variational space to arbitrary vector fields allowed us to split the vector-valued problem into a set of coupled scalar-valued problems for each component. Established solution procedures for such problems, as \SFEM\ and \DI, are adapted to this situation. Numerical experiments on the canonical unit sphere and surfaces with non-constant curvature established the consistency of all introduced methods.

The experiments further showed the tight interplay of topology, geometry, and dynamics. In all experiments the defect localization is related to the Gaussian curvature of the surface, $+1$ defects are found at extrema of the Gaussian curvature, while $-1$ defects are located at saddle points. We have further demonstrated the general possibility to reduce the overall energy by introducing additional defects and thus establishing non-trivial realizations of the Poincar\'e-Hopf theorem as energy minima. The proposed methods allow to further investigate this interplay. Here the effect of $\Kn$ as well as the impact of intrinsic and extrinsic contributions should be analyzed.

The introduced models and methods should also be complemented by more rigorous theoretical works on the convergence of the thin film limit. In analogy to scalar-valued problems an extension to evolving surfaces seems feasible. Beyond the mentioned fundamental issues, the model and methods are ready to be applied in the field of passive and active soft matter and surface bound, non-equilibrium physics comprising orientational order. Examples are passive \cite{Vitelli2006,LopezLeon2011,koning2013} and active \cite{Menzel2013} liquid crystals and polar fluids \cite{Ahmadi2006,Bois2011,Kruse2004} in thin shells, which are proposed models for a cell cortex \cite{Ramaswamy2016}.

Although the polar model, described by the Frank-Oseen energy and the introduced dynamic equations, already shows a variety of interesting effects, a nematic model will have additional features. Therefore, the Q-tensor Landau-de Gennes models should be focused on. With similar ideas of incorporating a tangentiality penalization a weak Q-tensor model on a surface could be derived and analyzed.
 
{\bf Acknowledgements}
This work is partially supported by the German Research Foundation through grant Vo889/18. We further acknowledge computing resources provided at JSC under grant HR06.

\appendix
\section{Thin film limit of penalized Frank-Oseen energy}
\label{sec:Limit}

Considering a thin shell \( \Omega_{\delta} = \surf\times[-\delta/2,\delta/2]\) around the surface \( \surf \) with thickness \( \delta \),
the local coordinates \( \lu \) and \( \lv \) of the surface immersion \( \Xb \) and an additional coordinate \( \lxi \),
which acts along the surface normal \(\surfNormal\),
lead to a thin shell parametrization
\( \XXb:U_{\delta} \rightarrow \R^{3} \) for the parameter domain $U_\delta:= U\times[-\delta/2,\delta/2]$, with $\XXb$ defined by
\begin{align}
  \XXb(\lu,\lv,\lxi) = \Xb(\lu,\lv) + \lxi\surfNormal(\lu,\lv) \formPeriod
\end{align}
The thickness \( \delta \) is sufficiently small to guarantee the injectivity of the pushforward, see \cite{Napoli2012b}.

For a better readability, we denote indices which mark all three components \( \{\lu,\lv,\lxi\} \) by capital letters.
The indices for the surface components \( \{\lu,\lv\} \) are denoted by small letters.
The metric tensor \( \extendDomain{\gb} \) of the thin shell is given by its components \( \gext_{IJ} = \partial_{I}\XXb \cdot \partial_{J}\XXb \), \ie,
\begin{align}
\begin{aligned}
  \gext_{ij} = g_{ij} - 2\lxi\shapeOperator_{ij} + \landauxin{2}_{ij} = g_{ij} + \landauxi_{ij} \formComma \quad
  \gext_{\lxi\lxi} = 1 \text{ and } \gext_{i\lxi} = \gext_{\lxi i} = 0 \formPeriod
\end{aligned}
\end{align}
The pure formal indices on $\landau$ extend the asymptotic polynomial behavior to tensor context and preserve summation conventions.
Hence, for the Christoffel symbols
\( \chext{I}{J}{K} = \frac{1}{2}\gext^{KL}\left( \partial_{I}\gext_{JL} + \partial_{J}\gext_{IL} - \partial_{L}\gext_{IJ} \right) \),
we obtain
\begin{align}
\begin{aligned}
  \chext{i}{j}{k} &= \ch{i}{j}{k} + \landauxi_{ij}^{k}, \quad
  \chext{i}{j}{\lxi} = \shapeOperator_{ij} + \landauxi_{ij}, \quad
  \chext{i}{\lxi}{k} = \chext{\lxi}{i}{k} = -\tensorfs{\shapeOperator}{i}{k} + \landauxi_{i}^{k}, \\
  \chext{\lxi}{\lxi}{K} &= \chext{I}{\lxi}{\lxi} = \chext{\lxi}{I}{\lxi} = 0 \formPeriod
\end{aligned}
\end{align}
We can approximate the square root of the determinant \( |\extendDomain{\gb}| \) on \( \surf \) by
\( \sqrt{|\extendDomain{\gb}|} = \sqrt{\gext_{\lxi\lxi}|\gb|} + \landauxi = \left( 1 + \landauxi \right)\sqrt{|\gb|} \).
Therefore, the volume element becomes
\begin{align}\label{eq:thinshellvolumeelement}
  \dV = \sqrt{|\extendDomain{\gb}|}d\lxi\wedge d\lu\wedge d\lv = \left( 1 + \landauxi \right)d\lxi \wedge \dS \formPeriod
\end{align}
The 3-tensor, with the same qualities as the volume element, is the Levi-Civita tensor
\begin{align}
  \lcext_{IJK} = \dV\left( \partial_{I}\XXb, \partial_{J}\XXb, \partial_{K}\XXb \right)
              = \sqrt{|\extendDomain{\gb}|}\varepsilon_{IJK}
              = \sqrt{|\gb|}\varepsilon_{IJK} + \landauxi_{IJK}\formComma
\end{align}
with the common Levi-Civita symbols \( \varepsilon_{IJK}\in\{-1,0,1\} \).
With the Levi-Civita tensor \( \tensor{\lc} \) on the surface, defined by
\( \lc_{ij} = \dS\left( \partial_{i}\Xb, \partial_{j}\Xb  \right) = \sqrt{|\gb|}\varepsilon_{ij} \),
and the fact, that all non-vanishing components of the Levi-Civita tensor \( \tensor{\lcext} \) in the thin shell have exactly one \( \lxi \)-index,
we obtain
\begin{align}\label{eq:lcexttolc}
  \lcext_{\lxi ij} = - \lcext_{i \lxi j} = \lcext_{ij \lxi} = \lc_{ij} + \landauxi_{ij} \formPeriod
\end{align}

For a better distinction, we use a semicolon in the thin shell and a straight line on the surface to mark the components of the covariant derivative,
\ie, for the vector fields \( \extendDomain{\pb}\in C^{1}\left( \Omega_{\delta}, \Tangent\Omega_{\delta} \right) \) and \( \pb\in C^{1}\left( \surf, \Tangent\surf \right) \),
we write
\begin{align}
  \tensorsf{\extendDomain{\p}}{I}{;J} &= \partial_{J}\extendDomain{\p}^{I} + \chext{J}{K}{I}\extendDomain{\p}^{K} \text{ and} \\
  \tensorsf{p}{i}{|j}     &= \partial_{j}p^{i} + \ch{j}{k}{i}p^{k} \formPeriod
\end{align}
The contravariant derivatives are given by
\( \extendDomain{\p}^{I;J} = \gext^{JK}\tensorsf{\extendDomain{\p}}{I}{;K} \)
and
\( p^{i|j} = g^{jk}\tensorsf{p}{i}{|k} \).
Henceforward, we assume that \( \extendDomain{\pb}\in\Tangent\Omega_{\delta} \) is an extension of \( \pb \), \ie, \( \extendDomain{\pb}\big|_{\surf} = \pb \in \Tangent\surf \),
and \( \extendDomain{\pb} \) is parallel and length-preserving in direction of \( \surfNormal \), \ie,
\( \tensorsf{\extendDomain{\p}}{I}{;\lxi} = 0 \) as a consequence\footnote{The constraints on \( \extendDomain{\pb} \) need to be physically interpreted and discussed. Other assumptions or boundary conditions on the outer shell surface finally lead to different models.}.
Therefore, the Taylor approximation on the surface of the contravariant tangential components becomes
\begin{align}
\begin{aligned}
  \extendDomain{\p}^{i} &= p^{i} + \lxi\partial_{\lxi}\extendDomain{\p}^{i}\big|_{\surf} + \landauxin{2}^{i}
             = p^{i} +\lxi\left( \tensorsf{\extendDomain{\p}}{i}{;\lxi} - \chext{\lxi}{K}{i}\extendDomain{\p}^{K}\big|_{\surf} \right) + \landauxin{2}^{i} \\
            &= p^{i} +\lxi\tensorfs{\shapeOperator}{k}{i}p^{k} + \landauxin{2}^{i} \formPeriod
\end{aligned}
\end{align}
It holds \( \extendDomain{\p}^{\lxi} = 0 \), because \( \extendDomain{\p}^{\lxi}\big|_{\surf} = 0 \) and
\( \partial_{\lxi}\extendDomain{\p}^{\lxi} = \tensorsf{\extendDomain{\p}}{\lxi}{;\lxi} - \chext{\lxi}{K}{\lxi}\extendDomain{\p}^{K} = 0 \),
but nonetheless, we get non-vanishing covariant tangential derivatives
\begin{align}
  \tensorsf{\extendDomain{\p}}{\lxi}{;j} = \chext{j}{K}{\lxi}\extendDomain{\p}^{K} = \shapeOperator_{jk}p^{k} + \landauxi_{j} \formPeriod
\end{align}
All remaining covariant derivatives can be approximated by
\begin{align}
  \tensorsf{\extendDomain{\p}}{i}{;j} &= \partial_{j}\extendDomain{\p}^{i} + \chext{j}{K}{i}\extendDomain{\p}^{K}
                           = \partial_{j}\p^{i} + \ch{j}{k}{i}\p^{k} + \landauxi^{i}_{j}
                          = \tensorsf{p}{i}{|j} + \landauxi^{i}_{j}  \formPeriod
\end{align}

The divergence of a vector field is the trace of its covariant derivative reads
\begin{align}\label{eq:divapprox}
  \nabla\cdot\extendDomain{\pb} &= \tensorsf{\extendDomain{\p}}{I}{;I}
                     = \tensorsf{\extendDomain{\p}}{i}{;i}
                     =  \tensorsf{p}{i}{|i} + \landauxi
                     = \Div\pb + \landauxi \formPeriod
\end{align}

The covariant curl of a vector field can be obtained by a double contraction of the Levi-Civita tensor and the contravariant derivative, \ie,
\begin{align}
 \left[ \nabla\times\extendDomain{\pb} \right]_{I} = -\lcext_{IJK}\extendDomain{\p}^{J;K} \formPeriod
\end{align}
With \eqref{eq:lcexttolc}, the \( \lxi \)-component of the curl can be approximated by
\begin{align}\label{eq:rotapproxnor}
  \left[ \nabla\times\extendDomain{\pb} \right]_{\lxi}
        &= -\lc_{jk}\gext^{kL} \tensorsf{\extendDomain{\p}}{j}{;L} + \landauxi
         = -\lc_{jk}g^{kl} \tensorsf{p}{j}{|l} + \landauxi
         = \Rot\pb + \landauxi
\end{align}
and the covariant tangential components by
\begin{align}\label{eq:rotapproxtan}
\begin{aligned}
  \left[ \nabla\times\extendDomain{\pb} \right]_{i}
      &= -\left( \lcext_{ij\lxi}\extendDomain{\p}^{j;\lxi}  +  \lcext_{i\lxi j}\extendDomain{\p}^{\lxi; j}\right)
       = \lc_{ij}\gext^{jK}\tensorsf{\extendDomain{\p}}{\lxi}{;K} + \landauxi_{i} \\
      &= \lc_{ij}\tensorsf{\shapeOperator}{j}{l}p^{l}  + \landauxi_{i}
       = -[*(\shapeOperator\pb)^{\flat}]_{i} + \landauxi_{i} \formComma
\end{aligned}
\end{align}
where we use, that for a every \( \qb\in\Tangent\surf \)
\begin{align}
  *\qb^{\flat} = \mathbf{i}_{\qb}(\!\dS) = \sqrt{|\gb|}\left( -q^{\lu}d\lv + q^{\lv}d\lu \right) = -\tensor{\lc}\qb
\end{align}
is valid on \( \surf \), see \cite{Marsden1988}.
The Hodge star operator is length-preserving and the metric \( \tensor{\gext} \) induces the common norm in the thin shell,
therefore it holds
\begin{align*}
  \left\| \nabla\times\extendDomain{\pb} \right\|^{2}_{\Omega_{\delta}}
        &= \left\| -*(\shapeOperator\pb)^{\flat} \right\|^{2}_{\surf} + \gext^{\lxi\lxi}\left(  \Rot\pb \right)^{2} + \landauxi
         = \left\| \shapeOperator\pb \right\|^{2}_{\surf} + \left(  \Rot\pb \right)^{2} + \landauxi \formPeriod
\end{align*}

Finally, with \( \left\| \extendDomain{\pb} \right\|^{2}_{\Omega_{\delta}} = \left\| \pb \right\|^{2}_{\surf} + \landauxi \),
\eqref{eq:thinshellvolumeelement}, \eqref{eq:divapprox}, \eqref{eq:rotapproxnor}, and \eqref{eq:rotapproxtan},
we can approximate the penalized Frank-Oseen energy \eqref{eq:FOEPenalized} in the thin shell \( \Omega_{\delta} \) by
\begin{align*}
  & \F{\Kn}\left[ \extendDomain{\pb}, \Omega_{\delta} \right] \\
      &= \int_{\surf} \int_{-\delta/2}^{\delta/2} \frac{\K}{2} \left( \left( \Div\pb \right)^{2} + \left( \Rot\pb \right)^{2} + \left\| \shapeOperator\pb \right\|^{2}_{\surf} \right)
                 + \frac{\Kn}{4}\left( \left\| \pb \right\|^{2}_{\surf} - 1 \right)^{2}  + \landauxi \; d\lxi\wedge\dS \\
      &=\delta\left( \E{\pb} + \landau(\delta) \right)
\end{align*}
for \( \extendDomain{\pb}\in\Hdr{\Omega_{\delta}}{\Tangent \Omega_{\delta}} \) and \( \pb\in\Hdr{\surf}{\Tangent \surf} \).

  \section{Integral Theorems}\label{sec:IntegralTheorems}
    The exterior derivative \( \exd \) is the \( L^{2} \)-adjoint of \( (-* \exd *) \). This allows to obtain some frequently used
    integral identities for the tangential vector field \( \pb=\alphab^{\sharp}:\surf\to\Tangent\surf \)
    on a closed surface \( \surf \) and
    also for its \( \R^{3} \) extension \( \pExt:\surf\to\R^3 \), with $ \pb = \ProjectSurf\pExt$.
    We get
    \begin{align*}
      -\int_{\surf} \Scalarprod{\Grad f , \pExt} \dS
              &= -\int_{\surf} \Scalarprod{\Grad f , \pb} \dS
              = -\int_{\surf} \Scalarprod{\exd f , \alphab} \dS\\
              &= \int_{\surf} f * \exd * \alphab \dS
             = \int_{\surf} f \Div \pb \dS \notag\\
             &= \int_{\surf} f \Div (\ProjectSurf\pExt) \dS
             = \int_{\surf} f \Div\pExt - \meanCurvature\left( \pExt \cdot \surfNormal \right) \dS\notag
    \end{align*}
    and
    \begin{align*}
      -\int_{\surf} \Scalarprod{\Rot f , \pExt} \dS
          &= -\int_{\surf} \Scalarprod{\Rot f , \pb} \dS
           = -\int_{\surf} \Scalarprod{*\exd f , \alphab} \dS\\
          &= \int_{\surf} \Scalarprod{\exd f , *\alphab} \dS
           = -\int_{\surf} f * \exd **\alphab \dS
           = \int_{\surf} f \Rot \pb \dS \notag \\
          &= \int_{\surf} f \Rot(\ProjectSurf\pExt) \dS
           = \int_{\surf} f \Rot \pExt \dS \notag \formPeriod
    \end{align*}
    Note that \( **\alpha = -\alpha \)
    and the inner product is invariant with respect to \( * \), \( \flat \), and \( \sharp \), applied to both arguments of the product simultaneously, see \cite{Marsden1988}.
    Hence, we obtain for the Laplace-DeRham operator
    \begin{align*}
      \int_{\surf} \Scalarprod{\laplaceDeRham\pb , \qExt} \dS
          &= \int_{\surf} \Scalarprod{\laplaceDeRham\pb , \qb} \dS
           = -\int_{\surf} \Scalarprod{\Grad\Div\pb , \qb} + \Scalarprod{\Rot\Rot\pb , \qb}\dS \\
          &= \int_{\surf} (\Div\pb)(\Div\qb) + (\Rot\pb)(\Rot\qb) \dS \notag\\
          &= \int_{\surf} \Div(\ProjectSurf\pExt)\Div(\ProjectSurf\qExt) + \Rot(\ProjectSurf\pExt)\Rot(\ProjectSurf\qExt) \notag\\
          &= \int_{\surf} (\Div\pExt - \meanCurvature\left( \pExt \cdot \surfNormal \right))(\Div\qExt - \meanCurvature\left( \qExt \cdot \surfNormal \right))
                            + (\Rot\pExt)(\Rot\qExt)  \notag
          \formPeriod
    \end{align*}

\section{Convergence study of the Laplace-deRham approximation}\label{sec:convApproxLaplaceDeRham}
To justify the approximation $\laplaceDeRham \pb \approx \laplaceDeRhamTilde \pExt + \Kt \surfNormal\left( \surfNormal \cdot \pExt \right) $ we set up a test case consisting of a vector-valued Helmholtz equation on an ellipsoidal surface $\ellipsoid$ (major axis: $1.0$, $0.5$, and $1.5$)
\begin{align}
 -\laplaceDeRham \pb + \pb  = -\laplaceDeRham \pb_s + \pb_s =: \vect{f} \quad \mbox{ on } \ellipsoid
\end{align}
with given analytical solution $\pb_{s} = \left[ -2y, 0.5 x, 0 \right]^T \in C( \ellipsoid, \Tangent\ellipsoid )$. We solve
\begin{align}
 -\laplaceDeRhamTilde \pExt + \pExt + \Kt \surfNormal \left( \surfNormal \cdot \pExt \right) = \vect{f} \quad \mbox{ on } \ellipsoid
\end{align}
using \SFEM\ on a conforming triangulation $\ellipsoid_h$ of $\ellipsoid$ with piecewise linear
Lagrange elements $\mathbb{V}_h(\ellipsoid_h) = \{ v_h \in C^0(\ellipsoid_h) \; : \; v_h|_T \in \mathbb{P}^1, \, \forall \, T \in \Fs \}$
as trial and test space for all components $\hat{\p}_{i}$. This leads to a sequence of linear discrete equations
\begin{multline}
 \int_{\ellipsoid_h} \NablaSurf \cdot \pExt \gDerivative_i \psi + \NablaSurf \cdot \left(\pExt \times \surfNormal \right) \NablaSurf \cdot \left( \vect{e}_i \psi \times \surfNormal \right) \dS \\
 +\int_{\ellipsoid_h}\hat{\p}_{i} \psi \dS +  \Kt \int_{\ellipsoid_h} \surfNormalI_{i}\left(\surfNormal \cdot\pExt\right) \psi \dS = \int_{\ellipsoid_h} f_i \psi \dS.
\end{multline}
To assemble and solve the resulting system we use the FEM-toolbox AMDiS \cite{AMDiS2007,AMDiS2015}.

\begin{figure}[!ht]
 \begin{center}
    \includegraphics{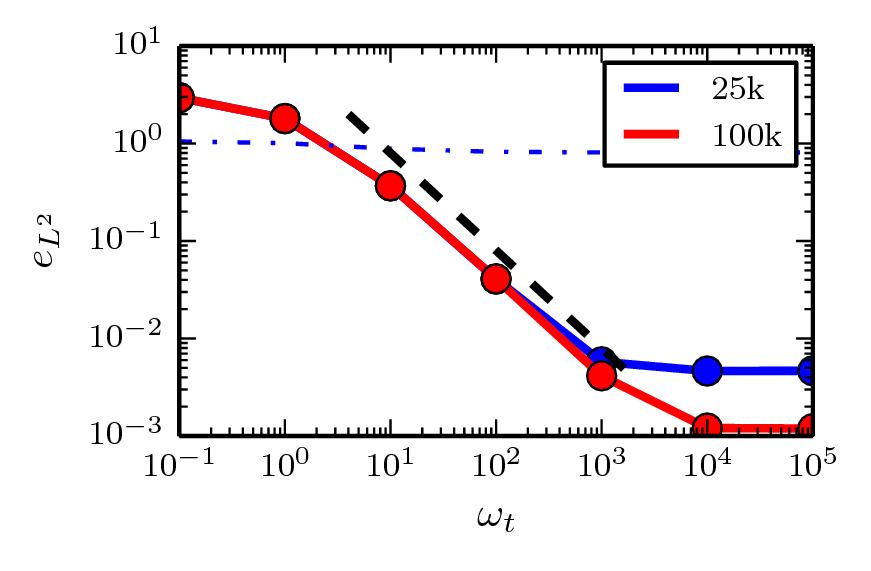}
  \caption{$L^2$-error for $\laplaceDeRhamTilde$ approximation (solid lines) for two well centered triangulations of $\ellipsoid$ with $25k$ and $100k$ vertices. The black dashed line indicates linear rate of convergence. The dash doted line shows the result for a component wise approximation of $\laplaceDeRham$ in \eqref{eq:LB}.}
  \label{fig:ConvergenceTestDeRhamL2}
 \end{center}
\end{figure}

\autoref{fig:ConvergenceTestDeRhamL2} shows the $L^2$-error $\epsilon_{L^2}(\pb) = \left( \int_{\ellipsoid} \sum_{i=1} (\pExt_i - \p_{s,i})^2 \dS  \right)^{1/2}$ vs  $\Kt$ and linear convergence, which is only limited by the mesh quality.\\

As a complementary result and to emphasize the delicate nature of the coupling between curvature and spatial derivatives, we also show in \autoref{fig:ConvergenceTestDeRhamL2} the $L^2$-error of a component wise approximation of $\laplaceDeRham$
\begin{align}
\label{eq:LB}
 \laplaceDeRham \pb \approx \sum_{i=i}^3 \NablaSurf \cdot \NablaSurf \hat{\p}_{i} \vect{e}_i + \Kt \surfNormal \left( \surfNormal \cdot \pExt \right).
\end{align}
As clearly visible in \autoref{fig:ConvergenceTestDeRhamL2}, this approximation fails for any values of $\Kt$ to reproduce the $\laplaceDeRham$ behavior on $\ellipsoid$.

\section{DEC: Notations and Details}
\label{sec:DECNotations}
\label{sec:DECLaplaceOperators}
\label{sec:DECLinOp}
 \subsection{Notations}

\begin{figure}
      \hspace{0.1\textwidth}
      \includegraphics[width=0.3\textwidth]{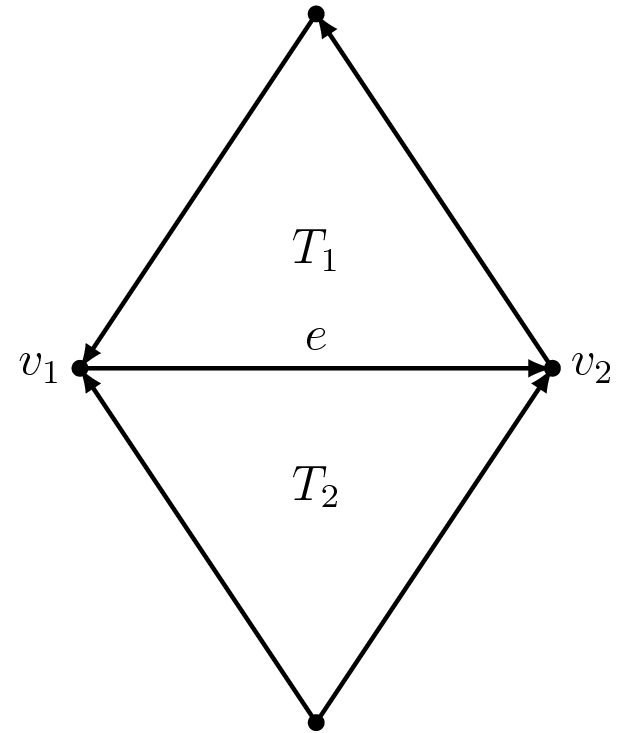}
      \hfill
      \includegraphics[width=0.35\textwidth]{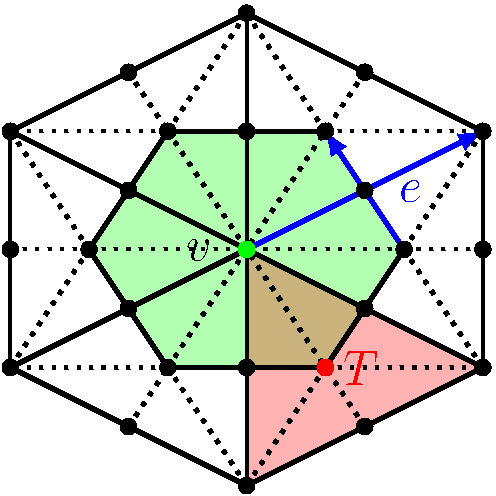}
      \hspace{0.1\textwidth}
      \caption{\textsc{Left:} This simple example mesh leads to
               \( s_{\face_{1},e} = +1 \), \( s_{\face_{2},e} = -1 \),
               \( s_{v_{1},e} = -1 \) and \( s_{v_{2},e} = +1 \).
               \textsc{Right:}
               the vertex \( v \) (green) and its Voronoi cell \( \star v \) (semi-transparent green);
               the edge \( e \) (blue) and its Voronoi edge \( \star e \) (blue);
               the face \( \face \) (semi-transparent red) and its Voronoi vertex (red).
               }
      \label{fig:orientationRelation}
    \end{figure}

    We often use the strict order relation \( \succ \) and \( \prec \) on simplices,
    where \( \succ \) is proverbial the ``contains'' relation,
    \ie, \( e \succ v \) means:
    the edge \( e \) contains the vertex \( v \).
    Correspondingly is \( \prec \) the ``part of'' relation,
    \ie, \( v \prec \face \) means:
    the vertex \( v \) is part of the face \( \face \).
    Hence, we can use this notation also for sums, like
    \( \sum_{\face\succ e} \), \ie, the sum over all faces \( \face \) containing edge \( e \),
    or \( \sum_{v\prec e} \), \ie, the sum over all vertices \( v \) being part of edge \( e \).
    Sometimes we need to determine this relation for edges more precisely with respect to the orientation.
    Therefore, a sign function is introduced,
    \begin{align}
      s_{\face,e} &:=
        \begin{cases}
          +1 & \text{if } e\prec \face \text{ and \( \face \) is on left side of \( e \)} \\
          -1 & \text{if } e\prec \face \text{ and \( \face \) is on right side of \( e \)} \\
           0 & e \nprec f\formComma
        \end{cases}\\
      s_{v,e} &:=
        \begin{cases}
          +1 & \text{if } v\prec e \text{ and \( e \) points to \( v \)} \\
          -1 & \text{if } v\prec e \text{ and \( e \) points away from \( v \)} \\
           0 & v \nprec e\formComma
        \end{cases}
    \end{align}
    to describe such relations between faces and edges, or vertices and edges, respectively.
    \autoref{fig:orientationRelation} gives a schematic picture.

    The property of primal mesh to be well-centered ensures the existence of a Voronoi mesh (dual mesh),
    which is also an orientable manifold-like simplicial complex, but not well-centered.

    The basis of the Voronoi mesh are not simplices, but chains of them.
    To identify these basic chains,
    we apply the (geometrical) star operator \( \star \) on the primal simplices,
    \ie, \( \star v\) is the Voronoi cell corresponding to the vertex \( v \)
    and inherits its orientation from the orientation of the polytope \( \left| \SC \right| \).
    \( \star v \) is, from a geometric point of view, the convex hull of circumcenters \( c(\face) \) of all triangles \( \face\succ v \).
    The Voronoi edge \( \star e \) of an edge \( e \) is a connection of the right face \( \face_{2}\succ e \) with the left face
    \( \face_{1}\succ e \) over the midpoint \( c(e) \).
    The Voronoi vertex \( \star\face \) of a face \( \face \) is simply its circumcenter \( c(\face) \)
    (see \autoref{fig:orientationRelation}).
    For greater details and a more mathematical discussion see, \eg, \cite{Hirani2003,Vanderzee2008}.

    The boundary operator \( \partial \) maps simplices (or chains of them) to the chain of simplices that describes its boundary,
    with respect to its orientation (see \cite{Hirani2003}),
    \eg, \( \partial(\star v)=-\sum_{e\succ v} s_{v,e}(\star e)\) (formal sum for chains)
    and \( \partial e = \sum_{v\prec e} s_{v,e} v \).

    The expression \( \left| \cdot \right| \) measures the volume of a simplex,
    \ie, \( \left| \face \right| \) the area of the face \( \face \),
    \( \left| e \right| \) the length of the edge \( e \)
    and the 0-dimensional volume \( \left| v \right| \) is set to be 1.
    Therefore, the volume is also defined for chains and the dual mesh, since the integral is a linear functional.

  \subsection{Laplace operators}
  With the Stokes theorem and the discrete Hodge operator defined in \cite{Hirani2003}
    we can develop a \DEC\ discretized Rot-Rot-Laplace for a discrete 1-form \( \alpha\in\Lambda_{h}^{1}(\SC) \) by
    \begin{align}
      \begin{aligned}
      \laplaceRotRot_h\alpha(e) &:=  \left(*\exd*\exd \alpha\right)(e)
                     = -\frac{\left| e \right|}{\left| \star e \right|} \left(\exd * \exd  \alpha\right)(\star e) \\
                    &=  -\frac{\left| e \right|}{\left| \star e \right|} \left(* \exd  \alpha\right)(\partial\star e)
                     = -\frac{\left| e \right|}{\left| \star e \right|} \sum_{\face\succ e} s_{\face,e} \left(* \exd \alpha\right) (\star
                     \face) \\
                    &= -\frac{\left| e \right|}{\left| \star e \right|} \sum_{\face\succ e} \frac{s_{\face,e}}{\left| \face \right|}
                                \left(\exd\alpha\right)(\face)
                     = -\frac{\left| e \right|}{\left| \star e \right|} \sum_{\face\succ e} \frac{s_{\face,e}}{\left| \face \right|}
                     \alpha(\partial \face)\\
                    &= -\frac{\left| e \right|}{\left| \star e \right|} \sum_{\face\succ e} \frac{s_{\face,e}}{\left| \face \right|}
                                \sum_{\tilde{e}\prec \face} s_{\face,\tilde{e}} \alpha(\tilde{e})
       \end{aligned}
    \end{align}
    and a \DEC\ discretized Grad-Div-Laplace by
    \begin{align}
      \begin{aligned}
      \laplaceGradDiv_h\alpha(e) &:= \left( \exd * \exd * \alpha\right) (e)
                      = \left( * \exd *\alpha \right) (\partial e) \\
                     &= \sum_{v\prec e} s_{v,e} \left( * \exd * \alpha \right)(v)
                      = \sum_{v\prec e} \frac{s_{v,e}}{\left| \star v \right|} \left( \exd * \alpha \right)(\star v)\\
                     &= \sum_{v\prec e} \frac{s_{v,e}}{\left| \star v \right|} \left( * \alpha \right)(\partial\star v)
                      = -\sum_{v\prec e} \frac{s_{v,e}}{\left| \star v \right|} \sum_{\tilde{e}\succ v}
                                s_{v,\tilde{e}} \left( * \alpha \right) (\star \tilde{e}) \\
                     &= -\sum_{v\prec e} \frac{s_{v,e}}{\left| \star v \right|} \sum_{\tilde{e}\succ v}
                                s_{v,\tilde{e}} \frac{\left| \star\tilde{e} \right|}{\left| \tilde{e} \right|}\alpha(\tilde{e}) \formPeriod
       \end{aligned}
    \end{align}
    Hence, we obtain the \DEC\ discretized Laplace-deRham operator by
    \begin{align*}
      \laplaceDeRham_h\alpha(e) = -\laplaceRotRot_h\alpha(e) - \laplaceGradDiv_h\alpha(e) \formPeriod
    \end{align*}

  \subsection{Conflate linear operators and its hodge dual to a PD-(1,1)-Tensor}

  For a linear operator \( \tensor{M}:\Tangent^{*}\surf\rightarrow\Tangent^{*}\surf \) point wise defined as a mixed co- and contravariant
  (1,1)-tensor with components \( \tensorfs{M}{i}{j} \), we discretize the 1-form \( \tensor{M}\alphab \) on an edge \( e\in\Es \) by definition \eqref{eq:discreteOneForm}
  and approximate the operator on the projected midpoint of the edge, \ie,
  \begin{align}\label{eq:MAlphaOnE}
    \left( \tensor{M}\alphab \right)_{h}(e)
      &= \int_{\pi(e)}\tensorfs{M}{i}{j}\alpha_{j}\ dx^{i}
      \approx \left[\tensor{M}(e)\right]_{ik}g^{kj}\int_{\pi(e)}\alpha_{j} d x^{i}\formComma
  \end{align}
  with \( \tensor{M}(e) := \tensor{M}|_{\pi(c(e))} \).
  With respect to an orthogonal basis \( \{\partial_{i}\xb,\partial_{j}\xb  \} \)
  with metric tensor \( \tensor{g}=g_{i}(dx^{i})^{2} \),
  we obtain for the 1-form \( \alphab=\alpha_{i}dx^{i} \) the Hodge dual
  \begin{align}
    *\alphab &= [*\alpha]_{1}dx^{1} + [*\alpha]_{2}dx^{2}
               = -\sqrt{\frac{g_{1}}{g_{2}}}\alpha_{2}dx^{1} + \sqrt{\frac{g_{2}}{g_{1}}}\alpha_{1}dx^{2} \formPeriod
  \end{align}
  Hence, we can replace the 1-forms beneath the integrals by
  \begin{align}\label{eq:differentialTrafos}
    \begin{bmatrix}
      \alpha_{1}dx^{1} & \alpha_{2}dx^{1} \\
      \alpha_{1}dx^{2} & \alpha_{2}dx^{2}
    \end{bmatrix}
      &=
    \begin{bmatrix}
      \alpha_{1}dx^{1} & -\sqrt{\frac{g_{2}}{g_{1}}}[*\alpha]_{1}dx^{1} \\
      \sqrt{\frac{g_{1}}{g_{2}}}[*\alpha]_{2}dx^{2} & \alpha_{2}dx^{2}
    \end{bmatrix} \formPeriod
  \end{align}
  Now, we use the basis \( \{\eb,\eb_{\star}\} \) defined in \autoref{sec:Dec} on the polytope \( |\SC| \)
  and the resulting metric \eqref{eq:PDMetric},
  \ie, \( g_{1}=|e|^{2} \) and \( g_{2} = |\star e|^{2} \).
  This leads to an approximation of \(  \left( \tensor{M}\alphab \right)_{h} \in \FormSpace_{h}(\SC) \)
  as a linear combination of \( \alpha_{h}, (*\alpha)_{h} \in \FormSpace_{h}(\SC)\),
  or rather, evaluated on an edge \( e\in\Es  \)
  \begin{align}\label{eq:Malpha}
    \left( \tensor{M}\alphab \right)_{h}(e)
      &\approx \frac{1}{|e|^{2}}M_{\eb,\eb}(e)\alpha_{h}(e) - \frac{1}{|e||\star e|}M_{\eb,\eb_{\star}}(*\alpha)_{h}(e)
  \end{align}
  and, in general, for \( \mathbf{v},\mathbf{w}\in\Span{\eb,\eb_{\star}} \) is
  \( M_{\mathbf{v},\mathbf{w}}(e) = \mathbf{v}\cdot\tensor{M}(e)\cdot\mathbf{w}
                                  = v^{i}\left[\tensor{M}(e)\right]_{ij}w^{j}\)
  the evaluation of the complete covariant tensor
  \( \tensor{M}(e) \) in direction \( \mathbf{v} \) and \( \mathbf{w} \).
  Note, if \( \tensor{M}\in\Tangent\surf\times\Tangent\surf\) is formulated in Euclidean \( \R^{3} \) coordinates,
  so that \( \tensor{M}(e)\in\R^{3 \times 3} \), there is no distinction between co- and contravariant components of \( \tensor{M}(e) \).
  Furthermore, if we use the approximation
  \( \left( *\tensor{M}\alphab \right)_{h}(e) \approx  -\frac{|e|}{|\star e|}\left( \tensor{M}\alphab \right)_{h}(\star e)\),
  we get with respect to \eqref{eq:MAlphaOnE} and \eqref{eq:differentialTrafos}
  \begin{align}\label{eq:starMalpha}
    \left( *\tensor{M}\alphab \right)_{h}(e)
        &\approx - \frac{1}{|e||\star e|}M_{\eb_{\star},\eb}\alpha_{h}(e) + \frac{1}{|\star e|^{2}}M_{\eb_{\star},\eb_{\star}}(*\alpha)_{h}(e)
        \formPeriod
  \end{align}
  Finally, we can summarize \eqref{eq:Malpha} and \eqref{eq:starMalpha} with the PD-1-form \( \alphav\in\FormSpace_{h}(\SC;\PDT^{*}\Es) \)
  on every edge \( e\in\Es \) to
  \begin{align}\label{eq:EndoApprox}
    \MPDt\cdot\alphav &:=
      \begin{bmatrix}
        \frac{1}{|e|^{2}}M_{\eb,\eb} & - \frac{1}{|e||\star e|}M_{\eb,\eb_{\star}} \\
        - \frac{1}{|e||\star e|}M_{\eb_{\star},\eb} & \frac{1}{|\star e|^{2}}M_{\eb_{\star},\eb_{\star}}
      \end{bmatrix}
      \cdot \alphav
      \approx
      \begin{bmatrix}
        \left( \tensor{M}\alphab \right)_{h} \\
        \left( *\tensor{M}\alphab \right)_{h}
      \end{bmatrix}
      \formComma
  \end{align}
  where the evaluation argument \( e \) is omitted for a better readability.

\newpage
\glsaddall
\printglossary[type=symbolslist]

\bibliographystyle{siam}
\bibliography{bibliography}
\end{document}